\documentclass[pra, aps, 10pt,twocolumn, floatfix]{revtex4-1}

\usepackage[titletoc,toc,title]{appendix}

\usepackage{amssymb}
\usepackage{amsmath}
\usepackage{amsfonts}
\usepackage{graphicx}
\usepackage{qtree}
\usepackage{float}

\usepackage{silence}
\usepackage{bm}

\usepackage{color}

\begin{document}

\title{Cumulants of Hawkes point processes}

\author{Stojan Jovanovi\'{c}}
\email{stojan.jovanovic@bcf.uni-freiburg.de}
\affiliation{Bernstein Center Freiburg \& Faculty of Biology, University of Freiburg, 79104 Freiburg im Breisgau, Germany}
\affiliation{KTH Royal Institute of Technology, 10691 Stockholm, Sweden}

\author{John Hertz}
\email{hertz@nbi.dk}
\affiliation{Institute for Neuroscience and Pharmacology and Niels Bohr Institute, University of Copenhagen, 2100 Copenhagen, Denmark}
\affiliation{NORDITA, KTH Royal Institute of Technology and Stockholm University, 10691 Stockholm, Sweden }

\author{Stefan Rotter}
\email{stefan.rotter@biologie.uni-freiburg.de}
\affiliation{Bernstein Center Freiburg \& Faculty of Biology, University of Freiburg, 79104 Freiburg im Breisgau, Germany }

\begin{abstract}
We derive explicit, closed-form expressions for the cumulant densities of a multivariate, self-exciting Hawkes point process, generalizing a result of Hawkes 
in his earlier work on the covariance density and Bartlett spectrum of such processes. To do this, we represent the Hawkes process in terms of a Poisson cluster
process and show how the cumulant density formulas can be derived by enumerating all possible "family trees", representing complex interactions between point 
events. We also consider the problem of computing the integrated cumulants, characterizing the average measure of correlated activity between events of different
types, and derive the relevant equations.
\end{abstract}

\maketitle
\thispagestyle{empty}

\section{Introduction}

The Hawkes point process was first introduced in \cite{Hawkes1} as a model of chain-reaction-like phenomena, in which the occurrence of an event increases the
likelihood of more such events happening in the future. This intrinsic ``self-exciting’’ property has made Hawkes processes appealing to a wide variety of researchers dealing 
with data exhibiting strong temporal clustering. Although it was originally used to model the dynamics of aftershocks that accompany strong earthquakes \cite{Ogata1988,Vere-Jones1970}, it has since found 
application to many other problems, including accretion disc formation \cite{Pechacek2008}, gene interactions \cite{Reynaud-Bouret2010}, social dynamics \cite{Mitchell2010}, insurance risk
\cite{Karabash2012,Zhu2013}, corporate default clustering \cite{Azizpour2008,AffinePointProcess}, market impact \cite{Bacry2014}, high-frequency financial data
\cite{Bauwens2009}, micro-structure noise \cite{Bacry2011}, crime \cite{Mohler2011}, generic properties of high-dimensional inverse problems in statistical mechanics \cite{AstronomyMastromatteo}, 
and dynamics of neural networks \cite{Pernice2011, PreOnaga}.  There is also recent theoretical work exploring generic mathematical properties of the process in its own respect 
\cite{PreSaichev1, PreSaichev2, PreBouchaud}.

As the areas of application of Hawkes models continue to grow, it becomes increasingly important to understand the probabilistic behavior of the process.  Unfortunately,
despite its ubiquity, the mathematical properties of the Hawkes process are still not fully known. In fact, the same dynamical characteristics that make it such a useful model
in practice are the ones that complicate formal analysis. Hawkes processes do not (except in some special cases, see \cite{DynamicContagion}) possess the Markov property,
making it impossible to study them using standard techniques. 

Recently, quite a few methods have been devised to circumvent this problem; there are now many well known results describing Hawkes process stability
\cite{Stability}, long-term behavior \cite{Karabash2012,Zhu2012} and large deviation properties \cite{Zhu2011}. Yet, since the early works of Hawkes himself on the covariance density and 
Bartlett spectrum \cite{Bartlett} of self-exciting processes \cite{Hawkes1,Hawkes2}, few have tried to further elucidate their statistical properties. In his work, Adamopoulos \cite{PGF}, 
for example, attempts to derive the probability generating functional of the Hawkes process, but manages only to represent it implicitly, as a solution of an intractable functional equation. 
Errais et al. \cite{AffinePointProcess}, using the elegant theory of affine jump processes, show that the moments of Hawkes processes can be computed by solving a system of non-linear
ODEs. Once again, however, explicit formulas turn out to be unobtainable by analytic means. Lastly, Saichev and Sornette \cite{Saichev2011,Saichev2011a}, using the alternative Poisson cluster 
representation of self-exciting processes, show that the moment generating function of the Hawkes process satisfies a transcendental equation which does not admit an explicit solution. 

Statistical behavior of, for example, Hawkes process moments and cumulants is of some importance in neuroscience, where the problem of quantifying levels of 
synchronization of action potentials has become very pertinent. It has been shown that nerve cells can be extremely sensitive to synchronous input from large groups
of neurons \cite{Rossant2011}. More precisely, a neuron's firing rate profile depends, to a large degree, on higher order correlations amongst the presynaptic spikes \cite{Kuhn2003}. 
Of course, which synchronous patterns are favored by the network is also determined by its connection structure. While the contribution of specific structural motifs to the emergence 
of pairwise correlations (i.e., two-spike patterns) has already been dissected \cite{Pernice2011}, no such result exists in the case of more complex patterns, 
stemming from correlations of higher order.

In this paper, we derive analytic formulas for the $n$th order cumulant densities of a linear, self-exciting Hawkes process with arbitrary interaction kernels, 
generalizing the result in \cite{Hawkes2}. Inspired by the approach of Saichev et al., we do this by utilizing the Poisson cluster process representation \cite{Cluster},
which simplifies calculations considerably. Furthermore, we show that the cumulant densities admit a natural and intuitive graphical representation in terms of the
branching structure of the underlying process and describe an algorithm that facilitates practical computation. Finally, we generalize the result in \cite{Pernice2011} by
showing that the integrated cumulant densities can be expressed in terms of formal sums of topological motifs of a graph, induced by specifying the physical interactions between
different types of point events.

\section{Preliminaries}

\subsection*{Basic definitions}

Consider a sequence $T = (T_{n})_{n\geq 1}$ of positive, random variables, representing times of random occurrences of a certain event. Alternatively, $T$ can be also thought of as as collection of 
random points on the positive half-line $\mathbb{R}^{+}$. By superposing all event times in the sequence, we obtain the point process $s = (s(t))_{t\geq0}$, formally defined by setting

\begin{equation}
 s(t):=\sum_{n\geq 1}\delta(t-T_{n})\textmd{,}
\end{equation}
where $\delta(t-T_{n})$ denotes the Dirac delta function, centered at the random point $T_{n}$. 

It is easy to see that the number of events occurring before time $t$ is given by

\begin{align}\label{count_process}
 N(t):= \int_{-\infty}^{t}s(u)du  = |\{n : T_{n}\leq t\}|\textmd{,}
\end{align}

The conditional probability, given the past activity, of a new event occurring in the interval $(t,t+dt)$ is given by the conditional rate function $(\lambda(t))_{t\geq0}$. More 
specifically, we have, up to first order \cite{CoxPoint}

\begin{equation}\label{rate_definition}
 P\{dN(t)=1 | \mathcal{H}_{t}\} = \lambda(t)dt\textmd{,}
\end{equation}
where $\mathcal{H}_{t}$ represents the history of the point process $s$ up to time $t$. Additionally, we assume that

\begin{equation}
  P\{dN(t)\geq 2 | \mathcal{H}_{t}\}=o(dt)\textmd{,}
\end{equation}
i.e. that the probability of two or more events arriving simultaneously is negligibly small. Intuitively, therefore, the conditional rate function represents the probability of a new event occurring in 
the infinitesimally near future, given the information about all events in the past. 

Furthermore, from our previous considerations it also follows that $dN(t)$ is (up to first order) a Bernoulli random variable and therefore,

 \begin{equation}\label{bernoulli}
  \langle dN(t) \rangle = P\{dN(t)=1\} = P\{\textmd{an event occurs at $t$}\}\textmd{.}
 \end{equation}

\subsection{The multivariate Hawkes process}
As was pointed out in \cite{Cluster}, the Hawkes process can be defined in two equivalent ways: either by specifying its conditional rate function
or as a \textit{Poisson cluster process}, generated by a certain branching structure.

\subsubsection{The conditional rate representation}
Following \cite{Hawkes1} and \cite{Hawkes2}, let us consider a $d$-dimensional point process $\mathbf{s} = (\mathbf{s}(t))_{t\geq0}$, with rate function 
$(\bm{\lambda}(t))_{t\geq0}$ \linebreak defined by

\begin{align}\label{rate}
 \bm{\lambda}(t) &:= \bm{\mu} + \int_{-\infty}^{t}\mathbf{G}(t-u)\cdot d\mathbf{N}(u)\\
 &\equiv\bm{\mu} + \int_{-\infty}^{t}\mathbf{G}(t-u)\cdot \mathbf{s}(u)du\textmd{,}
\end{align}
where $\bm{\mu}$ denotes the $d$-dimensional base rate vector with positive entries ($\mu^{i}>0$) and $\mathbf{G}(t)$ is an $d\times d$ matrix of non-negative, 
integrable functions $g^{ij}(t)$, with support on $\mathbb{R}^{+}$, called the interaction kernel. In principle, therefore, the rate $\bm{\lambda}(t)$ should 
always remain positive, but models for which the probability of negative values is sufficiently small may be useful approximations \cite{Volker}.

Rewriting equation (\ref{rate}) in terms of the components of the conditional rate function $\lambda^{i}(t)$, we find that, $\forall i$,

\begin{equation}\label{rate_comp}
 \lambda^{i}(t) = \mu^{i} + \sum_{j=1}^{d}\int_{-\infty}^{t}g^{ij}(t-u)dN^{j}(u)\textmd{.}
\end{equation}

From equations (\ref{rate_definition}) and (\ref{rate_comp}), we can now see that 

\begin{align}
 \frac{P\{dN^{i}(t) =1 |\mathcal{H}_{t}\}}{dt} &= \stackrel{\textmd{base rate}}{\overbrace{\mu^{i}}}\\ 
 &+ \underset{\textmd{influence of past $\mathcal{H}_{t}$}}{\underbrace{\sum_{j=1}^{d}\int_{-\infty}^{t}g^{ij}(t-u)dN^{j}(u)}}\textmd{,}
\end{align}
i.e. that the probability of an event of type $i$ occurring at time $t$ is simply the 
sum of a constant base rate and a convolution of the complete history of the process with the interaction kernel $\mathbf{G}(t)$, whose component $g^{ij}(t)$
describes the increase of the likelihood of type $i$ events at $t$, caused by a type $j$ event, occurring at $0$. Note that, in the special case of 
no interactions ($g^{ij}(t)\equiv0$), we recover the definition of a multivariate Poisson process with constant rate $\bm{\mu}$. In this case, however,
the (conditional) rate function is independent both of time and of the history $\mathcal{H}_{t}$.

\subsubsection{The cluster process representation}

\begin{figure}
 \begin{center}
  \includegraphics[width=0.5\textwidth]{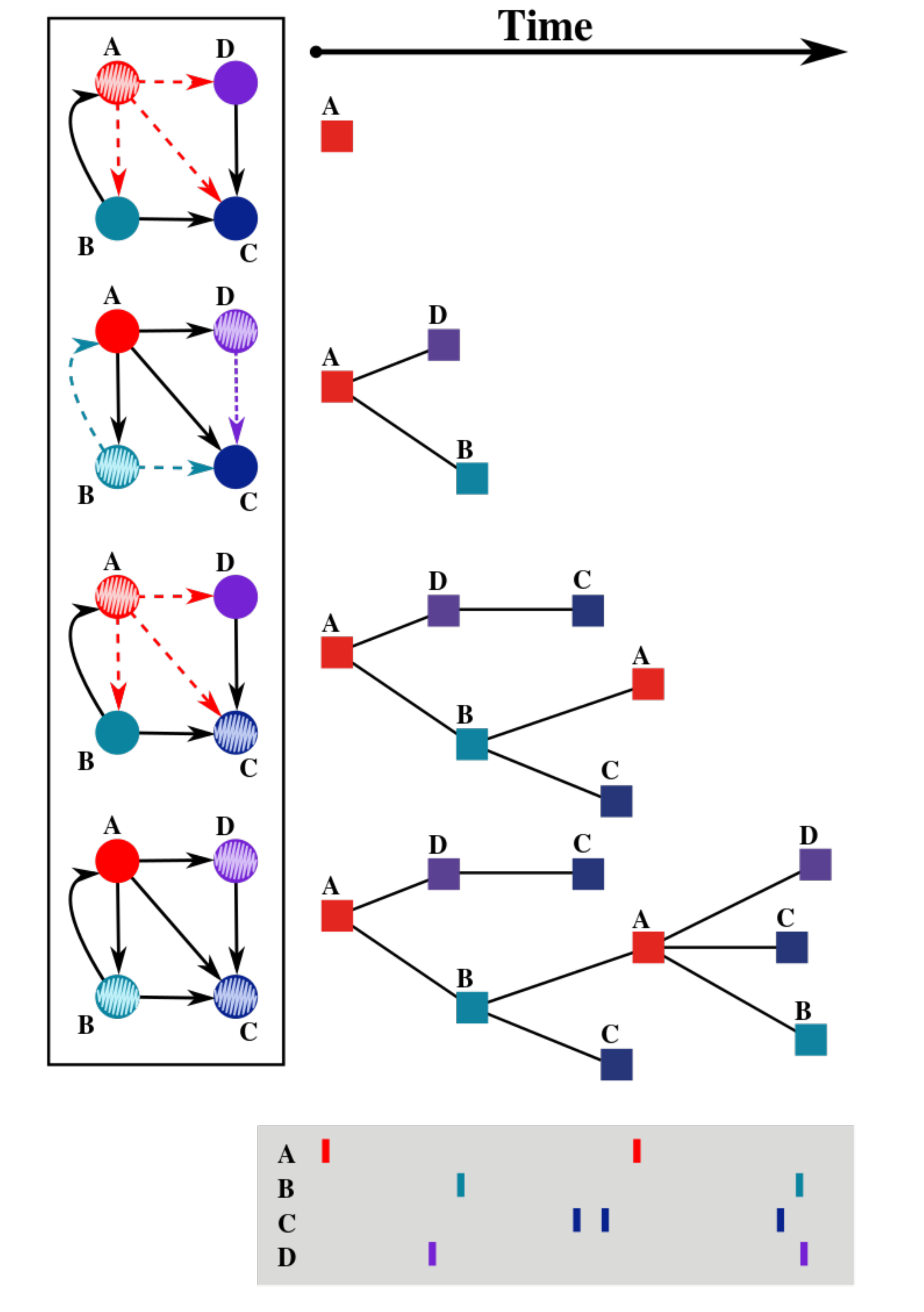}
  \caption{(Color online) Evolution of a Poisson cluster process on a network with $4$ nodes and $6$ directed links. \\Left column: Hatched nodes are active, and dashed links are transmitting a signal to all their respective
  neighbors.\\ Right column: A cluster, generated by the arrival of an immigrant from node $A$, evolves sequentially as new offspring are produced. \\Row $1$: A type $A$ immigrant arrives, providing the seed from which a new cluster will emerge. \\Row $2$: The type $A$ immigrant generates $1$st generation offspring of types $B$ and $D$, respectively.
  Together, these constitute the $1$st generation of events. \\Row $3$: The $1$st generation, type $D$, event generates a single offspring event of type $C$. The event of type $B$ creates two offspring of
  types $A$ and $C$, respectively. \\Row $4$: Finally, the $2$nd generation, type $A$, event is the only one to generate $3$ offspring events, of types $B$, $C$ and $D$, respectively. \\Bottom: $A$ time sequence 
  of events generated in this $4$ generation long evolution, here displayed by means of a "raster plot" which indicates the events generated by each node on the time axis.}
  \label{pedFigure}
 \end{center}
\end{figure}

Let us consider a Poisson cluster process $C$, which evolves in the following way (\cite{Bayes}, see also Figure \ref{pedFigure}) :

\begin{enumerate}
\item Let $I^{k}$ be a realization, on the interval $[0,T]$, of a homogeneous Poisson process with rate $\mu^{k}$. We will call points in $I^{k}$ \textit{immigrants of type $k$}.
\item For every $k$, each immigrant $x\in I^{k}$  generates a \textit{cluster} of points $C^{k}_{x}$. All such clusters are mutually independent.
\item The clusters $C^{k}_{x}$ are generated according to the \linebreak following branching structure:
			\begin{itemize}
			 \item Each cluster $C^{k}_{x}$ consists of generations of \linebreak offspring of all types of the immigrant $x$, which itself belongs to generation $0$.
			 \item Recursively, given the immigrant $x$ and the offspring of generation $1,2,\cdots,n$ of all types, every "child" $y$ of generation $n$ and 
			 type $j$, produces, $\forall i$, its own offspring of generation $n+1$ and type $i$ by generating a realization of an inhomogeneous Poisson process 
			 $O^{ij}$ with rate $\lambda(t):= g^{ij}(t - y)$. In other words, the probability of there being, at time $t$, a type $i$ offspring event of generation $n+1$, 
			 caused by a type $j$ event $y$ of generation $n$ is equal to $g^{ij}(t-y)dt$.
			\end{itemize}
\item The point process $C$ is equal to the superposition of all points in all generated clusters, i.e. 
			\begin{equation}
			 C = \sum_{k,x}C^{k}_{x}\textmd{.}
			\end{equation}
\end{enumerate}

For example, if  $C$ were used to model the dynamics of a spiking neuronal network, the immigrants $I^{k}$ would represent all the spikes of neuron $k$, caused by constant, external
input to the network, and the clusters $C^{k}_{x}$ all subsequent spikes, caused by action potential propagation through the network via synaptic connections.  

Having defined the cluster process $C$, it is then possible to show (see e.g. \cite{Cluster}) that by letting, $\forall i$ and $\forall t\geq0$,

\begin{equation}
 \lambda^{i}(t):= \lim_{\delta\rightarrow 0}\frac{1}{\delta}P\{\textmd{event of type $i$ in $[t, t+\delta)$}|\mathcal{H}_{t}\}\textmd{,}
\end{equation}
and assuming that the spectral radius $\rho(\mathbf{G})$ (i.e. the largest eigenvalue) of the integrated kernel matrix

\begin{equation}
 \mathbf{G}:=\int_{-\infty}^{+\infty}\mathbf{G}(t)dt
\end{equation}
is strictly less than $1$, then $\lambda^{i}(t)$ must be equal to the conditional rate function in equation (\ref{rate_comp}). Furthermore, 
the corresponding point process will also be stationary. In other words, we will have

\begin{equation}\label{stationary_rate}
 \frac{\langle \mathbf{N}(t)\rangle}{t} =  \langle \bm{\lambda}(t)\rangle = \bm{\lambda} = (\mathbf{I}-\mathbf{G})^{-1}\bm{\mu}\textmd{,}
\end{equation}
where $\mathbf{I}$ denotes the $d\times d$ identity matrix. In what follows, we will always assume that we are working with a stationary version of a Hawkes process.
More specifically, we will assume that $\rho(\mathbf{G})<1$ and, consequently, that the matrix $(\mathbf{I}-\mathbf{G})^{-1} $ can be expanded in terms of powers of the
integrated kernel matrix $\mathbf{G}$, i.e. 

\begin{equation}\label{expansion}
 (\mathbf{I}-\mathbf{G})^{-1} = \sum_{n=0}^{+\infty}\mathbf{G}^{n}\textmd{.}
\end{equation}

Note that the matrix $\mathbf{G}$ has a very useful interpretation (which follows from the definition of the Poisson cluster process) - its component 
$g^{ij}$ represents the average total number of events of type $i$ in the second generation, caused by a first generation, type $j$ event. Thus, the 
components of the $n$th matrix power $\mathbf{G}^{n}$ equal the average total number of type $i$ offspring within $n$ subsequent generations, of a first generation, 
type $j$ event.

From our previous considerations, it now follows that by requiring that $\rho(\mathbf{G})<1$ (or, equivalently, that the series (\ref{expansion}) converges) we, 
in fact, assume that each event of a given type produces only finitely many events of any other type, after an infinite number of generations.

\section{The Hawkes process cumulant density}

Consider now an arbitrary $n$-dimensional random vector $\mathbf{X} =(X_{1},\cdots, X_{n})\equiv X_{\bar{n}}$, where we used the symbol $\bar{n}$ to denote the set $\{1,\cdots,n\}$ The 
cumulant of order $n$, denoted by $k(X_{\bar{n}})$, is a general measure of statistical dependence of the components of $\mathbf{X}$. It is defined, combinatorially, as (see \cite{lukacs}, 
page 27.)
 
 \begin{equation}\label{cumulant}
 k(X_{\bar{n}}) = \sum_{\pi}(|\pi|-1)!(-1)^{|\pi|-1}\prod_{B\in\pi}\left\langle X_{B}\right\rangle\textmd{,}
\end{equation}
where the sum goes over all partitions $\pi$ of the set $\{1,\cdots,n\}$, $|\cdot|$ denotes the number of blocks of a given partition, and

\begin{equation}
 \langle X_{B} \rangle = \left\langle\prod_{i\in B}X_{i}\right\rangle\textmd{.}
\end{equation}

A dual formula, expressing moments in terms of cumulants, reads

\begin{equation}\label{moments_cumulant_formula}
 \langle X_{\bar{n}}\rangle = \sum_{\pi}\prod_{B\in\pi}k(X_{B})\textmd{,}
\end{equation}
where $k(X_{B})$ denotes the cumulant of those components of $\mathbf{X}$, whose indices are in $B$.

The cumulant $k(X_{\bar{n}})$ is a natural generalization, to higher dimensions, of the covariance $cov(X_{1},X_{2})$ of two variables. 

Indeed, if we set $n=2$, $\mathbf{X} = (X_{1},X_{2})$ and apply formula (\ref{cumulant}), we obtain

\begin{align}
 k(X_{1},X_{2}) &= (1-1)!(-1)^{1-1}\langle X_{1} X_{2}\rangle\nonumber\\ 
 &+ (2-1)!(-1)^{2-1} \langle X_{1}\rangle\langle X_{2}\rangle\textmd{,}
\end{align}
as $\pi_{1}=\{\{1,2\}\}$ and $\pi_{2}=\{\{1\},\{2\}\}$ are the only partitions of the set $\{1,2\}$. Also, obviously, $|\pi_{1}|=1$ and $|\pi_{2}|=2$. Thus,

\begin{equation}\label{covariance_cumulant}
  k(X_{1},X_{2})=\langle X_{1} X_{2}\rangle -  \langle X_{1}\rangle\langle X_{2}\rangle=cov(X_{1},X_{2})\textmd{.}
\end{equation}
 
For a given time vector $\mathbf{t} = (t_{1},\cdots, t_{n})$ and multi-index $\mathbf{i} = (i_{1},\cdots,i_{n})$, we now define the the 
$\mathbf{n}$th order cumulant density of the Hawkes process, denoted by $k^{\mathbf{i}}(\mathbf{t})$, by letting

\begin{equation}
 k^{\mathbf{i}}({\mathbf{t}}):= \frac{k(dN^{i_{1}}(t_{1}),\cdots,dN^{i_{n}}(t_{n}))}{d\mathbf{t}}\textmd{,}
\end{equation}
 where we used $d\mathbf{t}$ to denote the differential $dt_{1}\cdots dt_{n}$.
 
As in the general case, the cumulant density $k^{\mathbf{i}}(\mathbf{t})$ is used to quantify the mutual dependence of random events of types 
$(i_{1},\cdots,i_{n})$ at times $(t_{1},\cdots,t_{n})$. 

For example, from equations (\ref{covariance_cumulant}) and (\ref{bernoulli}), we can see that, for $\mathbf{i}=(1,2)$, $\mathbf{t}=(t_{1},t_{2})$ and $d\mathbf{t}=dt_{1}dt_{2}$,
 
 \begin{eqnarray}
  k^{\mathbf{i}}(\mathbf{t})d\mathbf{t}&&=P\{\textmd{type $1$ event at $t_{1}$, type $2$ event at $t_{2}$}\}\nonumber\\
  &&\hspace*{-25pt}-P\{\textmd{type $1$ event at $t_{1}$}\}P\{\textmd{type $2$ event at $t_{2}$}\}\textmd{.}
 \end{eqnarray}

The formulas for the $n$th order cumulant density $k^{\mathbf{i}}(\mathbf{t})$, however, get more and more complicated with increasing $n$,
as the number of set partitions involved grows supra-exponentially. 
 
To illustrate this point, we set $n=3$, $\mathbf{i} = (1,2,3)$ and $\mathbf{t}=(t_{1},t_{2},t_{3})$. Then, from (\ref{cumulant}), we have
 
  \begin{align}
  k^{\mathbf{i}}(\mathbf{t})d\mathbf{t} &= \langle dN^{1}(t_{1})dN^{2}(t_{2})dN^{3}(t_{3})\rangle\nonumber\\ 
  &- \langle dN^{1}(t_{1})dN^{2}(t_{2})\rangle\langle dN^{3}(t_{3})\rangle\nonumber\\
  &- \langle dN^{1}(t_{1})dN^{3}(t_{3})\rangle\langle dN^{2}(t_{2})\rangle\nonumber\\
  &- \langle dN^{2}(t_{2})dN^{3}(t_{3})\rangle\langle dN^{1}(t_{1})\rangle\nonumber\\
  &+2\langle dN^{1}(t_{1})\rangle\langle dN^{2}(t_{2})\rangle\langle dN^{3}(t_{3})\rangle\textmd{.}
 \end{align}

To alleviate the problem of increasing complexity, we use the cluster process representation to come up with a useful and intuitive expression for the density 
$k^{\mathbf{i}}(\mathbf{t})$ in terms of the cluster process's branching structure. 

First off, note that the only way that events $(t_{1},\cdots,t_{n})$ (of types $(i_{1},\cdots,i_{n})$) 
can be statistically dependent is if they all belong to the same cluster, i.e. if they are all offspring (possibly of different generations) of a single original immigrant.

More specifically, we can show that (see Appendix \ref{proof}), for every multi-index $\mathbf{i} = (i_{1},\cdots,i_{n})$ and every vector $\mathbf{t}=(t_{1},\cdots,t_{n})$,
 
 \begin{equation}\label{cumulant_probability_formula}
  k^{\mathbf{i}}(\mathbf{t})d\mathbf{t} = P\left\{E^{\mathbf{i}}_{\mathbf{t}}\cap C^{\mathbf{i}}_{\mathbf{t}}\right\}\textmd{,}
 \end{equation}
 where 
 
 \begin{align}
  &E^{\mathbf{i}}_{\mathbf{t}} = \{\textmd{$\forall k$, there is a type $i_{k}$ event at time $t_{k}$}\}\textmd{,}\\
  &C^{\mathbf{i}}_{\mathbf{t}} = \{ \exists\textmd{ cluster } C\textmd{ such that, }\forall k\textmd{, }t_{k}\in C\}\textmd{.}
 \end{align}

This result now provides us with a practical way of \linebreak computing $k^{\mathbf{i}}(\mathbf{t})d\mathbf{t}$. 

\begin{figure}
 \includegraphics[trim = 60mm 185mm 0mm 43mm]{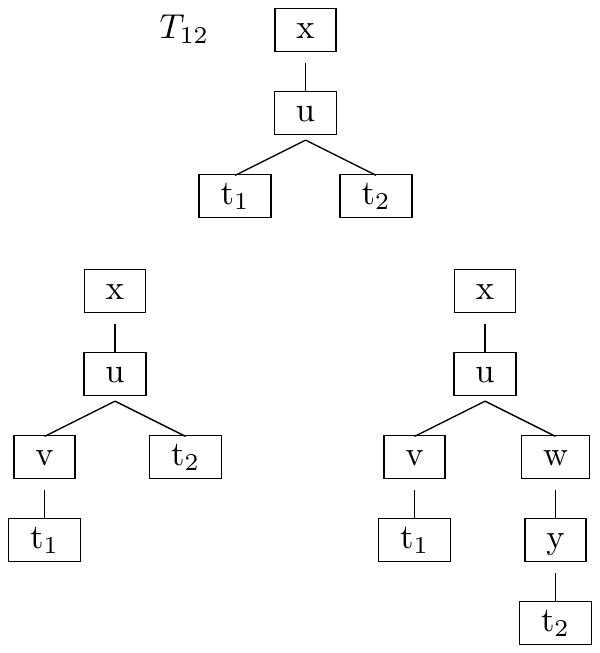}
 \caption{A schematic representation of all possible family trees, containing events $t_{1}$ and $t_{2}$. Two concrete examples of family trees, mapping to the same general schema are 
 given in the second row. The root, $x$, denotes the original immigrant and $u$ the branching point in the family
 tree of immigrant $x$, leading to the appearance of offspring $t_{1}$ and $t_{2}$. Note that each link connecting two nodes can be, in theory, any number of generations long.}
\label{2-trees}
\end{figure}

For example, in the case when $n=2$, we have that $k^{ij}(t_{1},t_{2})dt_{1}dt_{2}$ is equal to the probability of there being a type $i$ event at 
time $t_{1}$, a type $j$ event at time $t_{2}$, and that both of these events are descendant from a common immigrant. Therefore, in order to compute the $2$nd
order cumulant density, we need to sum up the probabilities of all possible "family trees" which contain events $t_{1}$ and $t_{2}$ (see Figure \ref{2-trees}).

% \begin{figure}
% \begin{center}
% \qtreecenterfalse
% \qtreeshowframes
% $T_{12}$ \Tree [.x [.u  t_{1} t_{2} !{\qbalance} ] ]
% \end{center}
% \caption{A schematic representation of all possible family trees, containing events $t_{1}$ and $t_{2}$. The root, $x$, denotes the original immigrant and $u$ the branching point in the family
% tree of immigrant $x$, leading to the appearance of offspring $t_{1}$ and $t_{2}$. Note that each link connecting two nodes can be, in theory, any number of generations long.}
% \label{2-trees}
% \end{figure}

In order to formalize this computation, we define 

\begin{equation}
R^{ij}_{t}:=\frac{P\{\textmd{type $j$ event at $0$ causes type $i$ event at $t$}\}}{dt}.
\end{equation}

Then,

\begin{equation}
 R^{ij}_{t} = \left[\sum_{n\geq0}\mathbf{G}^{\star n}(t)\right]_{ij}\textmd{,}
\end{equation}
where $[\cdot]_{ij}$ extracts component $(i,j)$ of a given matrix, and $\mathbf{G}^{\star n}(t)$ denotes the $\mathbf{n}$th convolution power of the 
interaction kernel $\mathbf{G}(t)$, defined recursively by

\begin{align*}
 \mathbf{G}^{\star 0}(t)&=\mathbf{I}\delta(t)\textmd{,}\\
 \mathbf{G}^{\star n}(t)&= \int_{-\infty}^{t}\mathbf{G}^{\star (n-1)}(t-s)\cdot \mathbf{G}(s)ds\textmd{.}
\end{align*}

Indeed, if we define $p^{ij}_{n}(t)$ to be equal to the probability that an event of type $j$ at $0$, after $n$ generations, 
causes a type $i$ event at $t$, we have
\begin{align}
 &\frac{p^{ij}_{0}(t)}{dt} = \delta_{ij}\delta(t)=[\mathbf{I}\delta(t)]_{ij}\textmd{,}\\ 
 &\frac{p^{ij}_{1}(t)}{dt} = [\mathbf{G}(t)]_{ij}\textmd{,}\\
 &\frac{p^{ij}_{2}(t)}{dt} = \sum_{k=1}^{d}\int_{-\infty}^{t}[\mathbf{G}(t-s)]_{ik}[\mathbf{G}(s)]_{kj}ds = [\mathbf{G}^{\star 2}(t)]_{ij}\textmd{,}
\end{align}
and therefore, by induction,

\begin{align}
 R^{ij}_{t} &= [I\delta(t)+\mathbf{G}(t)+\mathbf{G}^{\star 2}(t)+\mathbf{G}^{\star 3}(t)+\cdots]_{ij}\\
 &= \left[\sum_{n\geq0}\mathbf{G}^{\star n}(t)\right]_{ij}\textmd{.}
 \end{align}

Furthermore, noting that $P\{$type $k$ immigrant at $x\}$ is, by construction, equal to $\mu^{k} dx$, we obtain the probability 
of an immigrant (arriving at any point in time) generating an event of type $m$ at time $u$. It equals

\begin{equation}
 \sum_{k=1}^{d}\int_{\mathbb{R}}\mu^{k}R^{mk}_{u-x}dx = \sum_{k=1}^{d}[(\mathbf{I}-\mathbf{G})^{-1}]_{mk}\mu^{k} =\lambda^{m}\textmd{,}
\end{equation}
i.e. it is the $m$th component of the stationary rate vector $\bm{\lambda}$ in (\ref{stationary_rate}), where the first equality in the previous equation
follows from

\begin{align}
 \int_{\mathbb{R}}R^{mk}_{u-x}dx &= \sum_{n\geq 0}\int_{\mathbb{R}}[\mathbf{G}^{\star n}(u-x)]_{mk}dx \\
 &=\sum_{n\geq 0}[\mathbf{G}^{n}]_{mk} = [(\mathbf{I}-\mathbf{G})^{-1}]_{mk}\textmd{.} 
\end{align}

Computing the probability of the family tree $T_{12}$ in Figure \ref{2-trees} is now straightforward; recalling the definition of $R^{ij}_{t}$ and taking into account
our previous considerations, we get

\begin{align}
 &k^{ij}(t_{1},t_{2})=\frac{P(T_{12})}{dt_{1}dt_{2}}=\sum_{m=1}^{d}\lambda^{m}\int_{\mathbb{R}}R^{im}_{t_{1}-u}R^{jm}_{t_{2}-u}du\textmd{,}
\end{align}
recovering a classical and well known result on the \linebreak covariance density of the Hawkes process (see \cite{Hawkes2}).

A big advantage of our approach, however, is that it can be used to compute cumulant densities of orders greater than $2$. 

\begin{figure}
 \includegraphics[trim = 60mm 185mm 0mm 44mm, width = 0.8\textwidth]{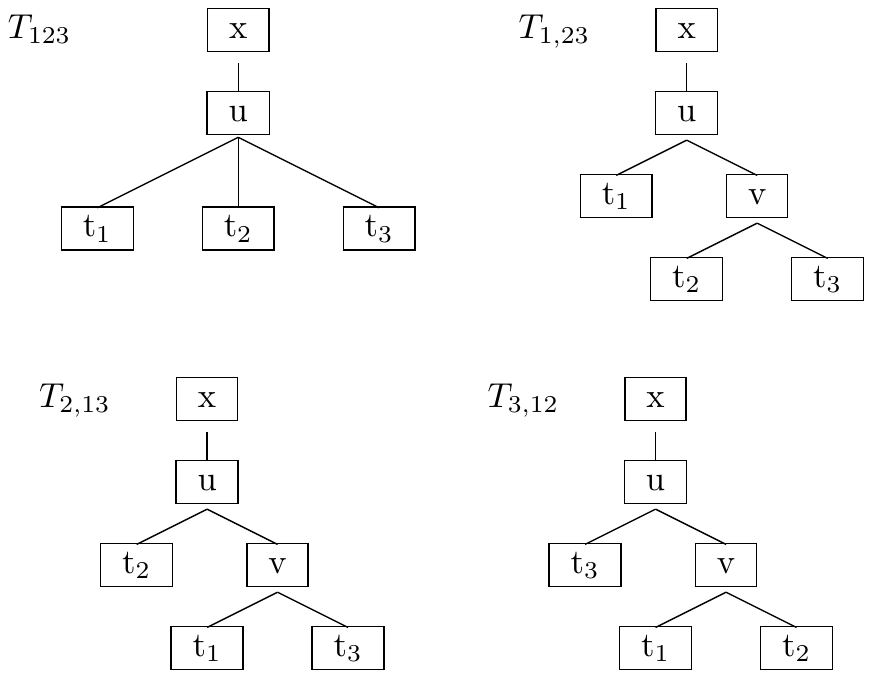}
 \caption{A schematic representation of all possible family trees, containing events $t_{1}$, $t_{2}$ and $t_{3}$. Once again, the root $x$ denotes the original 
 immigrant, nodes $u$ and $v$ represent the branching points in the family tree of immigrant $x$ and each link between two nodes can be any number of generations long.}
 \label{3-trees}
\end{figure}

For example, in order to compute the $3$rd order density $k^{ijk}(t_{1},t_{2},t_{3})$ we start, as in the $2$-dimensional case, by enumerating all possible 
family trees with leaves $t_{1}$, $t_{2}$ and $t_{3}$. In this case, however, there are in total $4$ different possibilities (see Figure \ref{3-trees}).

We can now proceed in much the same way as before, summing up the probabilities of all possible trees in order to derive the desired formula. We define $\mathbf{t}=(t_{1},t_{2},t_{3})$,
$d\mathbf{t}=dt_{1}dt_{2}dt_{3}$ and

\begin{equation}
 \Psi^{ij}_{t}=R^{ij}_{t}-\delta_{ij}\delta(t)=\left[\sum_{n\geq1}\mathbf{G}^{\star n}(t)\right]_{ij}\textmd{,}
\end{equation}
finally obtaining

\begin{align}\label{3-density}
 &k^{ijk}(\mathbf{t}) = \frac{P(T_{1,23})}{d\mathbf{t}}+\frac{P(T_{2,13})}{d\mathbf{t}} +\frac{P(T_{3,12})}{d\mathbf{t}}+\frac{P(T_{123})}{d\mathbf{t}}\nonumber\\
 &= \sum_{m,n=1}^{d}\lambda^{n}\int_{\mathbb{R}}R^{in}_{t_{1}-u}\left(\int_{\mathbb{R}}R^{jm}_{t_{2}-v}R^{km}_{t_{3}-v}\Psi^{mn}_{v-u}dv\right)du\nonumber\\
 &+ \sum_{m,n=1}^{d}\lambda^{n}\int_{\mathbb{R}}R^{jn}_{t_{2}-u}\left(\int_{\mathbb{R}}R^{im}_{t_{1}-v}R^{km}_{t_{3}-v}\Psi^{mn}_{v-u}dv\right)du\nonumber\\
 &+ \sum_{m,n=1}^{d}\lambda^{n}\int_{\mathbb{R}}R^{kn}_{t_{3}-u}\left(\int_{\mathbb{R}}R^{im}_{t_{1}-v}R^{jm}_{t_{2}-v}\Psi^{mn}_{v-u}dv\right)du\nonumber\\
 &+ \sum_{m=1}^{d}\lambda^{m}\int_{\mathbb{R}}R^{im}_{t_{1}-u}R^{jm}_{t_{2}-u}R^{km}_{t_{3}-u}du\textmd{.}
\end{align}

It is important to point out that equation (\ref{3-density}) can be derived in a different, albeit a more tedious way using martingale theory arguments, generalizing the derivation of Bacry et al. 
in \cite{bacry2012} for the second order cumulant density.

The newly introduced function $\Psi^{ij}_{t}$ corresponds to the probability of a type $j$ event at $0$ generating a type $i$ event at $t$, after 
at least one generation.

The appearance of such a term in the above equations is a consequence of the fact that, for instance, contracting the link
between nodes $u$ and $v$ in tree $T_{1,23}$ to a point turns it into $T_{123}$, which is already accounted for. Thus, in order to avoid counting certain 
configurations twice, we must introduce a "stiff" link between the two internal nodes $u$ and $v$ in trees $T_{1,23}$, $T_{2,13}$ and $T_{3,12}$. 

% \begin{figure}
% \begin{center}
% \qtreecenterfalse
% \qtreeshowframes
% $T_{123}$ \Tree [.x [.u  t_{1} t_{2} t_{3} !{\qbalance} ] ]
% \hskip 0.5in
% $T_{1,23}$ \Tree [.x [.u t_{1} [.v t_{2} t_{3} !{\qbalance} ] ] ]
% \hskip 0.5in
% \vskip 0.3in
% $T_{2,13}$ \Tree [.x [.u t_{2} [.v t_{1} t_{3} !{\qbalance} ] ] ]
% \hskip 0.5in
% $T_{3,12}$ \Tree [.x [.u t_{3} [.v t_{1} t_{2} !{\qbalance} ] ] ]
% \end{center}
% \caption{A schematic representation of all possible family trees, containing events $t_{1}$, $t_{2}$ and $t_{3}$. Once again, the root $x$ denotes the original 
% immigrant, nodes $u$ and $v$ represent the branching points in the family tree of immigrant $x$ and each link between two nodes can be any number of generations long.}
% \label{3-trees}
% \end{figure}

By generalizing the above considerations, it is possible to construct a general procedure for computing the $n$th order cumulant density 
$k^{i_{1}\cdots i_{n}}(t_{1},\cdots,t_{n})$.

\begin{enumerate}
 \item For a given $n\geq 2$, generate all possible rooted trees $T$ with $n$ leaves.
 \item Label the leaves of $T$ with ordered pairs $(i_{k},t_{k})$, in arbitrary order. Label the internal nodes (including the root) of $T$ arbitrarily.
 \item For every tree $T$, construct an integral term $I_{T}$, according the the following pseudo-algorithm :
	\begin{enumerate}
	  \item Set $I_{T} \leftarrow 1$;
	  \item For every edge in T, connecting a node $v$ of type $j_{v}$ to a leaf $t_{k}$ of type $i_{k}$ : 
		\begin{equation*}
		 I_{T} \leftarrow I_{T}\cdot R^{i_{k}j_{v}}_{t_{k}-v}dv\textmd{;}
		\end{equation*}

          \item For every edge in T, connecting an internal node $u$ of type $j_{u}$ to another internal node $v$ of type $j_{v}$ :
		\begin{equation*}
		  I_{T} \leftarrow I_{T}\cdot \Psi^{j_{v}j_{u}}_{v-u}du\textmd{;}
		\end{equation*}
          \item Let $x$ be the root of $T$. Set
		\begin{equation*}
		  I_{T} \leftarrow I_{T}\cdot \lambda^{j_{x}}\textmd{;}
		\end{equation*}
	  \item Integrate $I_{T}$ with respect to the variable $du$, for every internal node $u$.
	  \item Sum over all $j_{u}$ for every internal node $u$.
	  \item Sum over all $j_{x}$
	\end{enumerate}
 \item Add up all integral terms $I_{T}$ for every rooted tree $T$, generated in the first step, to obtain the $n$th order cumulant density.
\end{enumerate}

The principal difficulty of the above procedure lies its first step, i.e. in the enumeration of all topologically distinct rooted trees with $n$ labeled leaves. 
While there are known algorithms that can tackle this problem (see e.g. the classic text by Felsenstein \cite{Felsenstein}), the number of terms grows very quickly 
with increasing $n$ (see Figure \ref{number_of_terms}) and thus computing $k^{\mathbf{i}}(\mathbf{t})$ quickly becomes impractical.

\begin{figure}
\begin{center}
 \begin{tabular}{ |c|r| }
  \hline
  $n$ & Terms in $k^{\mathbf{i}}(\mathbf{t})$\\
  \hline
  $2$ & $1$ \\
  $3$ & $4$ \\
  $4$ & $26$ \\
  $5$ & $236$ \\
  $6$ & $2,752$ \\
  $7$ & $39,208$ \\
  $8$ & $660,302$ \\
  $9$ & $12,818,912$ \\
  $10$ & $282,137,824$ \\
% %   $11$ & $6,939,897,856$\\
% %   $12$ & $188,666,182,784$ \\
% %   $13$ & $5,617,349,020,544$ \\
% %   $14$ & $181,790,703,209,728$ \\
% %   $15$ & $6,353,726,042,486,272$\\
% %   $16$ & $238,513,970,965,257,728$\\
% %   $17$ & $9,571,020,586,419,012,608$\\
% %   $18$ & $408,837,905,660,444,010,496$\\
% %   $19$ & $18,522,305,410,364,986,906,624$\\
% %   $20$ & $887,094,711,304,119,347,388,416$\\
  \hline
\end{tabular}
\end{center}
\caption{Number of terms in $k^{\mathbf{i}}(\mathbf{t})$ for a given $n$ - from \cite{Felsenstein}}
\label{number_of_terms}
\end{figure}

\section{Integrated cumulants as sums of topological motifs}

 Let $k^{\mathbf{i}}(\mathbf{t})$ be, for a given time vector $\mathbf{t} = (t_{1},\cdots,t_{n})$ and multi-index $\mathbf{i} = (i_{1},\cdots,i_{n})$, 
 the $n$th order cumulant \linebreak density of a $d$-dimensional Hawkes process. We define the integrated cumulant of order $\mathbf{n}$, denoted simply by 
 $k^{\mathbf{i}}$, by setting
 
 \begin{align}
  k^{\mathbf{i}}:=\int_{\mathbb{R}^{n}_{+}}k^{\mathbf{i}}(\mathbf{t})d\mathbf{t}\textmd{.}
 \end{align}
 
 Note that $k^{\mathbf{i}}$ can be seen as the $n$-dimensional Laplace transform, "at zero", of $k^{\mathbf{i}}(\mathbf{t})$. Indeed, if we denote 
 
 \begin{equation}
  \mathcal{L}_{\bm{\omega}}(k^{\mathbf{i}}(\mathbf{t})) = \int_{\mathbb{R}^{n}_{+}}e^{-\bm{\omega}\cdot\mathbf{t}}k^{\mathbf{i}}(\mathbf{t})d\mathbf{t}\textmd{,}
 \end{equation}
 where $\bm{\omega} = (\omega_{1},\cdots,\omega_{n})\in\mathbb{C}^{n}$ and $\bm{\omega}\cdot\mathbf{t} = \sum_{i}\omega_{i}t_{i}$, we have, 
 clearly,

\begin{equation}
 k^{\mathbf{i}} = \mathcal{L}_{\mathbf{0}}(k^{\mathbf{i}}(\mathbf{t}))\textmd{.}
\end{equation}

Thus, if we define

\begin{equation}
\mathbf{R}_{t}:= 
 \begin{pmatrix}
  R^{11}_{t} & \cdots & R^{1d}_{t} \\
  \vdots & \ddots & \vdots \\
  R^{d1}_{t}& \cdots & R^{dd}_{t}
 \end{pmatrix}\textmd{ , }\bm{\Psi}_{t} = \mathbf{R}_{t}-\mathbf{I}\delta(t)\textmd{,}
\end{equation}
we can, by Laplace transforming the covariance density $k^{ij}(t_{1},t_{2})$, prove (see Appendix \ref{formulas_for_integrated_cumulants}) that,

\begin{equation}\label{integrated_covariance}
 k^{ij} = \sum_{m=1}^{d}\lambda^{m}[\mathbf{R}]_{im}[\mathbf{R}]_{jm}\textmd{,}
\end{equation}
where we set 

\begin{equation}
 \mathbf{R} = (\mathbf{I}-\mathbf{G})^{-1} = \mathcal{L}_{0}(\mathbf{R}_{t})\textmd{.}
\end{equation}

Expanding $\mathbf{R}$ in powers of $\mathbf{G}$, we get

\begin{equation}
 k^{ij} = \sum_{m=1}^{d}\sum_{k=0}^{+\infty}\sum_{l=0}^{+\infty}\lambda^{m}[\mathbf{G}^{k}]_{im}[\mathbf{G}^{l}]_{jm}\textmd{.}
\end{equation}

Interpreting now the matrix power $\mathbf{G}^{l}$ in the sense of graph theory, i.e. as a matrix whose component $(i,j)$ corresponds to the sum of 
lengths of all paths from node $j$ to node $i$ in exactly $l$ steps, we see that the integrated covariance density $k^{ij}$ can be equivalently represented as

\begin{equation}
 k^{ij} = \sum_{T\in\mathcal{T}^{m}_{ij}}w(T)\textmd{,}
\end{equation}
where the sum goes over the set $\mathcal{T}^{m}_{ij}$ of all rooted trees $T$ with root $m$, containing nodes $i, j$. Here, $w(T)$ \linebreak denotes the 
weight of tree $T$, defined as the product of weights of all edges, contained in $T$, times the weight of the root $m$, defined as being equal 
to $\lambda^{m}$.

The graph $H$ with adjacency matrix $\mathbf{G}$ can be thought of as follows. Each node $i\in\{1,\cdots,n\}$ in $H$ corresponds to a type of event in the 
underlying Hawkes process, and the existence of an edge $e_{ij}$ from $j$ to $i$ indicates the possibility of generating type $i$ events from those of
type $j$. Starting in node $j$, traversing the corresponding edge to reach node $i$ is equivalent to generating $g^{ij}$ type $i$ offspring of a type $j$ immigrant.
Therefore, each path through graph $H$ represents a specific "bloodline" of a type $m$
immigrant, while a tree $T\in\mathcal{T}^{m}_{ij}$ accounts for the possibility of the bloodline splitting somewhere along the way, concluding in, after a 
certain number of generations, in offspring of both types $i$ and $j$. The previous formula tells us that the sum of weights of all such trees is equal to the
integrated covariance $k^{ij}$.

Now, reasoning in much the same way as before we have, for $k^{ijk}$,

\begin{align}\label{integrated_thrid_cumulant}
k^{ijk} &= \sum_{m=1}^{d}\lambda_{m}[\mathbf{R}]_{im}[\mathbf{R}]_{jm}[\mathbf{R}]_{km}\\
&+ \sum_{m,n=1}^{d}\lambda_{n}[\mathbf{R}]_{im}[\mathbf{R}]_{jm}[\bm{\Psi}]_{mn}[\mathbf{R}]_{kn}\\
&+ \sum_{m,n=1}^{d}\lambda_{n}[\mathbf{R}]_{jm}[\mathbf{R}]_{km}[\bm{\Psi}]_{mn}[\mathbf{R}]_{in}\\
&+ \sum_{m,n=1}^{d}\lambda_{n}[\mathbf{R}]_{im}[\mathbf{R}]_{km}[\bm{\Psi}]_{mn}[\mathbf{R}]_{jn}\textmd{,}
\end{align}
where $\bm{\Psi} = (\mathbf{I}-\mathbf{G})^{-1}-\mathbf{I} = \mathcal{L}_{0}(\bm{\Psi}_{t})$. Once again, expanding $\mathbf{R}$ and $\bm{\Psi}$
in powers of $\mathbf{G}$ yields

\begin{equation}
 k^{ijk} = \sum_{T\in\mathcal{T}^{m}_{ijk}}w(T)\textmd{,}
\end{equation}
where $\mathcal{T}^{m}_{ijk}$ is the set of all rooted trees with root $m$, containing nodes $i,j,k$ and $w(\cdot)$ is the already defined weight function.

It is now easy to see that the general result is of the form

\begin{equation}
 k^{\mathbf{i}} = \sum_{T\in\mathcal{T}^{m}_{\mathbf{i}}}w(T)\textmd{,}
\end{equation}
where $\mathbf{i} = (i_{1},\cdots,i_{n})$ and $\mathcal{T}^{m}_{\mathbf{i}} = \mathcal{T}^{m}_{i_{1}\cdots i_{n}}$ is the set of all rooted trees with root $m$, containing
nodes $i_{1},\cdots,i_{n}$.

\section{Discussion}

In this paper we described the method for computing a class of statistics of linear Hawkes self-exciting point processes with arbitrary interaction kernels. 
By using the Poisson cluster process representation, we were able to obtain a general procedure for deriving formulas for $n$th order cumulant densities. 
Furthermore, we have shown there is a one-to-one correspondence between the integral terms, appearing in said densities, and all topologically distinct rooted 
trees with $n$ labeled leaves.

We also considered the problem of computing time-integrated cumulants and showed this can be done by simplifying the expressions for the corresponding cumulant
densities. Moreover, and not surprisingly, we demonstrated that integrated cumulants likewise admit a representation in terms of a formal sum of topological 
motifs, generalizing previous work on the topological expansion of the integrated covariance \cite{Pernice2011}.

The problem of quantifying higher-order correlations is of some importance in theoretical neuroscience. Indeed, it has long been suggested \cite{hebb1949, gerstein1989} that 
understanding the cooperative dynamics of populations of neurons would provide fundamental insight into the nature of neuronal computation. However, while direct experimental evidence for coordinated 
activity on the spike train level mostly relies on the correlations between pairs of nerve cells \cite{Gray1989, vaadia1995, riehle1997, bair2001, kohn2005},
it is becoming increasingly clear that such pairwise correlations cannot completely resolve the cooperative dynamics of neuronal populations \cite{martignon1995, bohte2000, Kuhn2003, Nature2010} and 
that higher-order cumulants need to be taken into account. 

One possible shortcoming of our work is the (supra-exponentially) increasing complexity of the closed-form expressions for the densities $k^{\mathbf{i}}(\mathbf{t})$
for higher values of $n$. This "explosion", however, is mostly due to combinatorial factors, that arise in many problems involving cumulants. As their definition
naturally involves objects such as set partitions, it seems to us that these sorts of issues would be quite difficult to avoid.

Another limitation of the present model is that it only allows for excitatory interactions - an arrival of an event at a given time can only \textit{increase}
the likelihood of future event, never decrease it. We hope, in the future, to be able to extend our analysis to include models in which there also exists a 
possibility of mutual inhibition between points of different types.

Further generalizations of our results might involve computing cumulants (and other important statistics) of a non-linear Hawkes processes (see e.g. 
\cite{Stability} for the definition), whose conditional rate function involves a non-linear transformation of equation (\ref{rate}), thus allowing for, for example,
multiplicative interaction between point events \cite{Cardanobile2009}. However, in this case, the resulting process no longer admits an immigrant-offspring representation, meaning
an alternative approach would be necessary.

\acknowledgments{Supported by the Erasmus Mundus Joint Doctoral programme EuroSPIN and the German Federal Ministry of Education and Research (BFNT - Freiburg*T\"{u}bingen, grant 01GQ0830). Stojan Jovanovi\'{c} 
acknowledges the hospitality of NORDITA. The authors would also like to thank the referees for many helpful comments and suggestions, Marcel Sauerbier for useful discussion and Gunnar Grah for making Figure \ref{pedFigure}.}

\begin{appendices}
 \section{Proof of equation (\ref{cumulant_probability_formula})}\label{proof}
 
 Let $\mathbf{t}$ be an arbitrary time vector $\mathbf{t} = (t_{1},\cdots, t_{n})$ and $\mathbf{i}$ an arbitrary multi-index $\mathbf{i} = (i_{1},\cdots,i_{n})$.
 
 From (\ref{bernoulli}), for every vector $\mathbf{t} = (t_{1},\cdots, t_{n})$ and multi-index $\mathbf{i} = (i_{1},\cdots,i_{n})$ we have that
 
 \begin{align}
  &\langle dN^{i_{1}}(t_{1})\cdots dN^{i_{n}}(t_{n})\rangle= P\{E^{\mathbf{i}}_{\mathbf{t}}\}\textmd{,}\\
  &P\{E^{\mathbf{i}}_{\mathbf{t}}\} = P\{\textmd{$\forall k$, there is a type $i_{k}$ event at $t_{k}$}\}\textmd{.}
 \end{align}
 
 Furthermore, it is clear that
 
 \begin{equation}\label{decomposition_of_average}
  P\{E^{\mathbf{i}}_{\mathbf{t}}\} = P\{E^{\mathbf{i}}_{\mathbf{t}}\cap C^{\mathbf{i}}_{\mathbf{t}}\}+P\{E^{\mathbf{i}}_{\mathbf{t}}\cap \bar{C}^{\mathbf{i}}_{\mathbf{t}}\}\textmd{,}
 \end{equation}
 where  $\bar{C}^{\mathbf{i}}_{\mathbf{t}}$ denotes the complement of the set

\begin{equation}
 C^{\mathbf{i}}_{\mathbf{t}} = \{ \exists\textmd{ cluster } C\textmd{ such that, }\forall k\textmd{, }t_{k}\in C\}\textmd{.}
\end{equation}

Indeed, events $\mathbf{t}$ of type $\mathbf{i}$ either are, or aren't all in some cluster $C$. We now proceed by induction in $n$. For $n=2$, we have
 
 \begin{equation}
  P\{E^{ij}_{t_{1}t_{2}}\} = P\{E^{ij}_{t_{1}t_{2}}\cap C^{ij}_{t_{1}t_{2}}\} + P\{E^{ij}_{t_{1}t_{2}}\cap \bar{C}^{ij}_{t_{1}t_{2}}\}\textmd{.}
 \end{equation}

 But, as the only way that two events are not in the same cluster is if they each belong to a different one; say, if $t_{1}\in C$ and $t_{2}\in D$,
 
 \begin{align}
 &P\{E^{ij}_{t_{1}t_{2}}\cap \bar{C}^{ij}_{t_{1}t_{2}}\}=\sum_{C,D}P\{E^{i}_{t_{1}}\cap E^{j}_{t_{2}}\cap C^{i}_{t_{1}}\cap D^{j}_{t_{2}}\}\\
&=\sum_{C,D}P\{E^{i}_{t_{1}}\cap C^{i}_{t_{1}}\}P\{E^{j}_{t_{2}}\cap D^{j}_{t_{2}}\}\\
&=P\{E^{i}_{t_{1}}\}P\{E^{j}_{t_{2}}\} = \langle dN^{i}_{t_{1}}\rangle\langle dN^{j}_{t_{2}}\rangle\textmd{,}
\end{align}
because of independence of different clusters $C$ and $D$. Thus,

\begin{align}
 P\{E^{ij}_{t_{1}t_{2}}\cap C^{ij}_{t_{1}t_{2}}\} &= P\{E^{ij}_{t_{1}t_{2}}\} - \langle dN^{i}(t_{1})\rangle\langle dN^{j}(t_{2})\rangle\\
 &=k^{ij}(t_{1},t_{2})dt_{1}dt_{2}\textmd{,}
\end{align}
proving that formula (\ref{cumulant_probability_formula}) is true for $n=2$.

Next, we assume that (\ref{cumulant_probability_formula}) is true for $n-1$ and prove that it then must also be true for $n$.
 
Consider the complementary set $\bar{C}^{\mathbf{i}}_{\mathbf{t}}$. If events $\mathbf{t}$ are not all in the same cluster, how could they be distributed? One 
possibility is that they are divided up between two different clusters, like in the previous case. In fact, they could potentially be distributed in $c$ different 
clusters, where $2\leq c \leq n$. Therefore,

\begin{equation}
 P\{E^{\mathbf{i}}_{\mathbf{t}}\cap \bar{C}^{\mathbf{i}}_{\mathbf{t}}\} = \sum_{c=2}^{n}\prod_{r=1}^{c}P\{E^{\mathbf{i}_{r}}_{\mathbf{t}_{r}}\cap B(r)^{\mathbf{i}_{r}}_{\mathbf{t}_{r}}\}\textmd{,}
\end{equation}
where the first sum goes over all possible numbers of different clusters that events $\mathbf{t}$ could be partitioned in, while $\mathbf{t}_{r}$ denotes the subset 
of $\mathbf{t}$ that belong to the $r$th cluster $B(r)$ (and $\mathbf{i}_{r}$ denotes their types).

Now, note that the previous equation is, in fact, a sum over all partitions $\pi$ of the set $\{1,\cdots,n\}$ with at least two blocks (i.e. $|\pi|>1$). Let us now 
fix one such partition $\pi = \{B(1),\cdots, B(|\pi|)\}$. As $|\pi|>1$, we must have, $\forall r$, $1\leq r \leq |\pi|$, that $|B(r)|<n$. But then, by the inductive 
assumption,

\begin{equation}
 P\{E^{\mathbf{i}_{r}}_{\mathbf{t}_{r}}\cap B(r)^{\mathbf{i}_{r}}_{\mathbf{t}_{r}}\} = k^{\mathbf{i}_{r}}(\mathbf{t}_{r})\textmd{,}
\end{equation}
and, therefore, 

\begin{equation}\label{different_clusters}
 P\{E^{\mathbf{i}}_{\mathbf{t}}\cap \bar{C}^{\mathbf{i}}_{\mathbf{t}}\} = \sum_{\pi : |\pi|>1}\prod_{B\in \pi}k^{\mathbf{i}_{r}}(\mathbf{t}_{r})d\mathbf{t}_{r}\textmd{.}
\end{equation}

Finally, from (\ref{decomposition_of_average}), (\ref{different_clusters}) and (\ref{moments_cumulant_formula}), we get 

\begin{equation}
 k^{\mathbf{i}}_{\mathbf{t}}d\mathbf{t} + \sum_{\pi : |\pi|>1}\prod_{B\in \pi}k^{\mathbf{i}_{r}}(\mathbf{t}_{r})d\mathbf{t}_{r} = P\{E^{\mathbf{i}}_{\mathbf{t}}\cap C^{\mathbf{i}}_{\mathbf{t}}\} + P\{E^{\mathbf{i}}_{\mathbf{t}}\cap \bar{C}^{\mathbf{i}}_{\mathbf{t}}\}\textmd{,}
\end{equation}
which completes the proof.

 \section{Formulas for integrated cumulants}\label{formulas_for_integrated_cumulants}
 
 Let $\mathcal{T}^{m}_{\mathbf{i}}$ be the set of all rooted trees with root $m$ and leaves $\mathbf{i} = (i_{1},\cdots, i_{n})$. Next, let $T\in\mathcal{T}^{m}_{\mathbf{i}}$
 and let $I_{T}$ be the corresponding integral term. In order to compute the Laplace transform $\mathcal{L}_{\bm{\omega}}(I_{T})$, we first consider the 
 leaves of $T$.
 
 Each leaf $i_{k}$ contributes a term $R^{i_{k}j_{v}}_{t_{k}-v}$, for some internal node $v$. For simplicity, let us assume that leaves $i_{1},\cdots,i_{s(v)}$
 all descend from a single internal node, which we denote $v$. Then, applying to $I_{T}$ the Laplace transform with respect to variables $t_{1},\cdots, t_{s(v)}$, we
 obtain
 
 \begin{equation}\label{leaf_terms}
  \prod_{l=1}^{s(v)}\mathcal{L}_{\omega_{l}}(R^{i_{l}j_{v}}_{t_{l}-v})e^{-v\sigma_{v}}\textmd{,}
 \end{equation}
 where we denote $\sigma_{v}=\sum_{l=1}^{s(v)}\omega_{l}$.

 Of course, in general the leaves $i_{1},\cdots,i_{n}$ are divided into several groups, according to which internal node they descend from. In that case, applying
 to each such group the Laplace transform in the already described way, yields several terms of type (\ref{leaf_terms}).
 
 Moving one level up in tree $T$, we are now in a situation in which several internal nodes, each with its own group of dependent leaves, all descend from a common
 node $u$, residing one level above them. We denote these internal nodes by $v_{1},\cdots, v_{r(u)}$. Each such internal node $v_{l}$ contributes to $I_{T}$ a term
 $\bm{\Psi}^{j_{v_{l}}j_{u}}_{v_{l}-u}$. Transforming the exponential term in (\ref{leaf_terms})
 
 \begin{equation}
  e^{-v_{l}\sigma_{v_{l}}} = e^{-(v_{l}-u)\sigma_{v_{l}}}e^{-u\sigma_{v_{l}}}\textmd{,}
 \end{equation}
 and multiplying with $\bm{\Psi}^{j_{v_{l}}j_{u}}_{v_{l}-u}$, we get 
 
 \begin{equation}\label{internal_node_terms}
  \prod_{l=1}^{r(u)}\mathcal{L}_{\sigma_{v_{l}}}(R^{j_{v_{l}}j_{u}}_{v_{l}-u})e^{-u\Sigma_{u}}\textmd{,}
 \end{equation}
 where $\Sigma_{u} = \sum_{l=1}^{r(u)}\sigma_{v_{l}}$.
 
 By induction, we can now see that this procedure must end after a finite number of steps (equal to the "depth" of tree $T$), at which point we are left with a 
 product of various terms of types (\ref{leaf_terms}) and (\ref{internal_node_terms}), integrated with respect to the position $x$ of the root $m$ (as this is the
 last node we reach by "climbing up" $T$). The exponential terms in this product can be combined to form
 
 \begin{equation}
  \int_{\mathbb{R}}e^{-x\sum_{i=1}^{n}\omega_{i}}dx = \delta(\omega_{1}+\cdots+\omega_{n})\textmd{,}
 \end{equation}
the integral representation of a Dirac delta function. 

By setting $\bm{\omega} = \mathbf{0}$, we now see that the formulas for $k^{\mathbf{i}}$ can be obtained from formulas for the cumulant densities by simply
"erasing" all the integral signs and replacing all the functional terms with their integrated counterparts.
 
\end{appendices}

% \bibliography{Bib/bibl}

\begin{thebibliography}{49}%
\makeatletter
\providecommand \@ifxundefined [1]{%
 \@ifx{#1\undefined}
}%
\providecommand \@ifnum [1]{%
 \ifnum #1\expandafter \@firstoftwo
 \else \expandafter \@secondoftwo
 \fi
}%
\providecommand \@ifx [1]{%
 \ifx #1\expandafter \@firstoftwo
 \else \expandafter \@secondoftwo
 \fi
}%
\providecommand \natexlab [1]{#1}%
\providecommand \enquote  [1]{``#1''}%
\providecommand \bibnamefont  [1]{#1}%
\providecommand \bibfnamefont [1]{#1}%
\providecommand \citenamefont [1]{#1}%
\providecommand \href@noop [0]{\@secondoftwo}%
\providecommand \href [0]{\begingroup \@sanitize@url \@href}%
\providecommand \@href[1]{\@@startlink{#1}\@@href}%
\providecommand \@@href[1]{\endgroup#1\@@endlink}%
\providecommand \@sanitize@url [0]{\catcode `\\12\catcode `\$12\catcode
  `\&12\catcode `\#12\catcode `\^12\catcode `\_12\catcode `\%12\relax}%
\providecommand \@@startlink[1]{}%
\providecommand \@@endlink[0]{}%
\providecommand \url  [0]{\begingroup\@sanitize@url \@url }%
\providecommand \@url [1]{\endgroup\@href {#1}{\urlprefix }}%
\providecommand \urlprefix  [0]{URL }%
\providecommand \Eprint [0]{\href }%
\providecommand \doibase [0]{http://dx.doi.org/}%
\providecommand \selectlanguage [0]{\@gobble}%
\providecommand \bibinfo  [0]{\@secondoftwo}%
\providecommand \bibfield  [0]{\@secondoftwo}%
\providecommand \translation [1]{[#1]}%
\providecommand \BibitemOpen [0]{}%
\providecommand \bibitemStop [0]{}%
\providecommand \bibitemNoStop [0]{.\EOS\space}%
\providecommand \EOS [0]{\spacefactor3000\relax}%
\providecommand \BibitemShut  [1]{\csname bibitem#1\endcsname}%
\let\auto@bib@innerbib\@empty
%</preamble>
\bibitem [{\citenamefont {Hawkes}(1971{\natexlab{a}})}]{Hawkes1}%
  \BibitemOpen
  \bibfield  {author} {\bibinfo {author} {\bibfnamefont {A.~G.}\ \bibnamefont
  {Hawkes}},\ }\href@noop {} {\bibfield  {journal} {\bibinfo  {journal}
  {Biometrika}\ }\textbf {\bibinfo {volume} {58}},\ \bibinfo {pages} {83}
  (\bibinfo {year} {1971}{\natexlab{a}})}\BibitemShut {NoStop}%
\bibitem [{\citenamefont {Ogata}(1988)}]{Ogata1988}%
  \BibitemOpen
  \bibfield  {author} {\bibinfo {author} {\bibfnamefont {Y.}~\bibnamefont
  {Ogata}},\ }\href {\doibase 10.1080/01621459.1988.10478560} {\bibfield
  {journal} {\bibinfo  {journal} {Journal of the American Statistical
  Association}\ }\textbf {\bibinfo {volume} {83}},\ \bibinfo {pages} {9}
  (\bibinfo {year} {1988})}\BibitemShut {NoStop}%
\bibitem [{\citenamefont {Vere-Jones}(1970)}]{Vere-Jones1970}%
  \BibitemOpen
  \bibfield  {author} {\bibinfo {author} {\bibfnamefont {D.}~\bibnamefont
  {Vere-Jones}},\ }\href@noop {} {\bibfield  {journal} {\bibinfo  {journal}
  {Journal of the Royal Statistical Society. Series B (Methodological)}\
  }\textbf {\bibinfo {volume} {32}},\ \bibinfo {pages} {pp. 1} (\bibinfo {year}
  {1970})}\BibitemShut {NoStop}%
\bibitem [{\citenamefont {Pech{\'a}cek}\ \emph {et~al.}(2008)\citenamefont
  {Pech{\'a}cek}, \citenamefont {Karas},\ and\ \citenamefont
  {Czerny}}]{Pechacek2008}%
  \BibitemOpen
  \bibfield  {author} {\bibinfo {author} {\bibfnamefont {T.}~\bibnamefont
  {Pech{\'a}cek}}, \bibinfo {author} {\bibfnamefont {V.}~\bibnamefont {Karas}},
  \ and\ \bibinfo {author} {\bibfnamefont {B.}~\bibnamefont {Czerny}},\
  }\href@noop {} {\bibfield  {journal} {\bibinfo  {journal} {Astronomy and
  Astrophysics}\ }\textbf {\bibinfo {volume} {487}},\ \bibinfo {pages} {815}
  (\bibinfo {year} {2008})}\BibitemShut {NoStop}%
\bibitem [{\citenamefont {Reynaud-Bouret}\ and\ \citenamefont
  {Schbath}(2010)}]{Reynaud-Bouret2010}%
  \BibitemOpen
  \bibfield  {author} {\bibinfo {author} {\bibfnamefont {P.}~\bibnamefont
  {Reynaud-Bouret}}\ and\ \bibinfo {author} {\bibfnamefont {S.}~\bibnamefont
  {Schbath}},\ }\href {\doibase 10.1214/10-AOS806} {\bibfield  {journal}
  {\bibinfo  {journal} {The Annals of Statistics}\ }\textbf {\bibinfo {volume}
  {38}},\ \bibinfo {pages} {2781} (\bibinfo {year} {2010})},\ \Eprint
  {http://arxiv.org/abs/arXiv:0903.2919v4} {arXiv:arXiv:0903.2919v4}
  \BibitemShut {NoStop}%
\bibitem [{\citenamefont {Mitchell}\ and\ \citenamefont
  {Cates}(2010)}]{Mitchell2010}%
  \BibitemOpen
  \bibfield  {author} {\bibinfo {author} {\bibfnamefont {L.}~\bibnamefont
  {Mitchell}}\ and\ \bibinfo {author} {\bibfnamefont {M.~E.}\ \bibnamefont
  {Cates}},\ }\href@noop {} {\bibfield  {journal} {\bibinfo  {journal} {Journal
  of Physics A: Mathematical and Theoretical}\ }\textbf {\bibinfo {volume}
  {43}},\ \bibinfo {pages} {045101} (\bibinfo {year} {2010})}\BibitemShut
  {NoStop}%
\bibitem [{\citenamefont {Karabash}\ and\ \citenamefont
  {Zhu}(2012)}]{Karabash2012}%
  \BibitemOpen
  \bibfield  {author} {\bibinfo {author} {\bibfnamefont {D.}~\bibnamefont
  {Karabash}}\ and\ \bibinfo {author} {\bibfnamefont {L.}~\bibnamefont {Zhu}},\
  }\href@noop {} {\bibfield  {journal} {\bibinfo  {journal} {arXiv preprint
  arXiv:1211.4039v2}\ ,\ \bibinfo {pages} {1}} (\bibinfo {year} {2012})},\
  \Eprint {http://arxiv.org/abs/arXiv:1211.4039v2} {arXiv:arXiv:1211.4039v2}
  \BibitemShut {NoStop}%
\bibitem [{\citenamefont {Zhu}(2013)}]{Zhu2013}%
  \BibitemOpen
  \bibfield  {author} {\bibinfo {author} {\bibfnamefont {L.}~\bibnamefont
  {Zhu}},\ }\href@noop {} {\bibfield  {journal} {\bibinfo  {journal}
  {Insurance: Mathematics and Economics}\ ,\ \bibinfo {pages} {1}} (\bibinfo
  {year} {2013})},\ \Eprint {http://arxiv.org/abs/arXiv:1304.1940v2}
  {arXiv:arXiv:1304.1940v2} \BibitemShut {NoStop}%
\bibitem [{\citenamefont {Azizpour}\ and\ \citenamefont
  {Giesecke}(2008)}]{Azizpour2008}%
  \BibitemOpen
  \bibfield  {author} {\bibinfo {author} {\bibfnamefont {S.}~\bibnamefont
  {Azizpour}}\ and\ \bibinfo {author} {\bibfnamefont {K.}~\bibnamefont
  {Giesecke}},\ }\href@noop {} {\bibfield  {journal} {\bibinfo  {journal}
  {Management Science}\ ,\ \bibinfo {pages} {1}} (\bibinfo {year}
  {2008})}\BibitemShut {NoStop}%
\bibitem [{\citenamefont {Errais}\ \emph {et~al.}(2010)\citenamefont {Errais},
  \citenamefont {Giesecke},\ and\ \citenamefont
  {Goldberg}}]{AffinePointProcess}%
  \BibitemOpen
  \bibfield  {author} {\bibinfo {author} {\bibfnamefont {E.}~\bibnamefont
  {Errais}}, \bibinfo {author} {\bibfnamefont {K.}~\bibnamefont {Giesecke}}, \
  and\ \bibinfo {author} {\bibfnamefont {L.~R.}\ \bibnamefont {Goldberg}},\
  }\href {\doibase 10.1137/090771272} {\bibfield  {journal} {\bibinfo
  {journal} {SIAM Journal on Financial Mathematics}\ }\textbf {\bibinfo
  {volume} {1}},\ \bibinfo {pages} {642} (\bibinfo {year} {2010})}\BibitemShut
  {NoStop}%
\bibitem [{\citenamefont {Bacry}\ and\ \citenamefont {Muzy}(2014)}]{Bacry2014}%
  \BibitemOpen
  \bibfield  {author} {\bibinfo {author} {\bibfnamefont {E.}~\bibnamefont
  {Bacry}}\ and\ \bibinfo {author} {\bibfnamefont {J.}~\bibnamefont {Muzy}},\
  }\href@noop {} {\bibfield  {journal} {\bibinfo  {journal} {arXiv preprint
  arXiv:1401.0903}\ ,\ \bibinfo {pages} {1}} (\bibinfo {year} {2014})},\
  \Eprint {http://arxiv.org/abs/arXiv:1401.0903v1} {arXiv:arXiv:1401.0903v1}
  \BibitemShut {NoStop}%
\bibitem [{\citenamefont {Bauwens}\ and\ \citenamefont
  {Hautsch}(2009)}]{Bauwens2009}%
  \BibitemOpen
  \bibfield  {author} {\bibinfo {author} {\bibfnamefont {L.}~\bibnamefont
  {Bauwens}}\ and\ \bibinfo {author} {\bibfnamefont {N.}~\bibnamefont
  {Hautsch}},\ }\href@noop {} {\emph {\bibinfo {title} {Modelling financial
  high frequency data using point processes}}}\ (\bibinfo  {publisher}
  {Springer},\ \bibinfo {year} {2009})\BibitemShut {NoStop}%
\bibitem [{\citenamefont {Bacry}\ \emph {et~al.}(2011)\citenamefont {Bacry},
  \citenamefont {Delattre}, \citenamefont {Marc},\ and\ \citenamefont
  {Muzy}}]{Bacry2011}%
  \BibitemOpen
  \bibfield  {author} {\bibinfo {author} {\bibfnamefont {E.}~\bibnamefont
  {Bacry}}, \bibinfo {author} {\bibfnamefont {S.}~\bibnamefont {Delattre}},
  \bibinfo {author} {\bibfnamefont {H.}~\bibnamefont {Marc}}, \ and\ \bibinfo
  {author} {\bibfnamefont {J.-F.}\ \bibnamefont {Muzy}},\ }in\ \href@noop {}
  {\emph {\bibinfo {booktitle} {Acoustics, Speech and Signal Processing
  (ICASSP), 2011 IEEE International Conference on}}}\ (\bibinfo {organization}
  {IEEE},\ \bibinfo {year} {2011})\ pp.\ \bibinfo {pages}
  {5740--5743}\BibitemShut {NoStop}%
\bibitem [{\citenamefont {Mohler}\ \emph {et~al.}(2011)\citenamefont {Mohler},
  \citenamefont {Short}, \citenamefont {Brantingham}, \citenamefont
  {Schoenberg},\ and\ \citenamefont {Tita}}]{Mohler2011}%
  \BibitemOpen
  \bibfield  {author} {\bibinfo {author} {\bibfnamefont {G.~O.}\ \bibnamefont
  {Mohler}}, \bibinfo {author} {\bibfnamefont {M.~B.}\ \bibnamefont {Short}},
  \bibinfo {author} {\bibfnamefont {P.~J.}\ \bibnamefont {Brantingham}},
  \bibinfo {author} {\bibfnamefont {F.~P.}\ \bibnamefont {Schoenberg}}, \ and\
  \bibinfo {author} {\bibfnamefont {G.~E.}\ \bibnamefont {Tita}},\ }\href
  {\doibase 10.1198/jasa.2011.ap09546} {\bibfield  {journal} {\bibinfo
  {journal} {Journal of the American Statistical Association}\ }\textbf
  {\bibinfo {volume} {106}},\ \bibinfo {pages} {100} (\bibinfo {year}
  {2011})}\BibitemShut {NoStop}%
\bibitem [{\citenamefont {Mastromatteo}\ and\ \citenamefont
  {Marsili}(2011)}]{AstronomyMastromatteo}%
  \BibitemOpen
  \bibfield  {author} {\bibinfo {author} {\bibfnamefont {I.}~\bibnamefont
  {Mastromatteo}}\ and\ \bibinfo {author} {\bibfnamefont {M.}~\bibnamefont
  {Marsili}},\ }\href {http://stacks.iop.org/1742-5468/2011/i=10/a=P10012}
  {\bibfield  {journal} {\bibinfo  {journal} {Journal of Statistical Mechanics:
  Theory and Experiment}\ }\textbf {\bibinfo {volume} {2011}},\ \bibinfo
  {pages} {P10012} (\bibinfo {year} {2011})}\BibitemShut {NoStop}%
\bibitem [{\citenamefont {Pernice}\ \emph
  {et~al.}(2011{\natexlab{a}})\citenamefont {Pernice}, \citenamefont {Staude},
  \citenamefont {Cardanobile},\ and\ \citenamefont {Rotter}}]{Pernice2011}%
  \BibitemOpen
  \bibfield  {author} {\bibinfo {author} {\bibfnamefont {V.}~\bibnamefont
  {Pernice}}, \bibinfo {author} {\bibfnamefont {B.}~\bibnamefont {Staude}},
  \bibinfo {author} {\bibfnamefont {S.}~\bibnamefont {Cardanobile}}, \ and\
  \bibinfo {author} {\bibfnamefont {S.}~\bibnamefont {Rotter}},\ }\href
  {\doibase 10.1371/journal.pcbi.1002059} {\bibfield  {journal} {\bibinfo
  {journal} {PLoS computational biology}\ }\textbf {\bibinfo {volume} {7}},\
  \bibinfo {pages} {e1002059} (\bibinfo {year}
  {2011}{\natexlab{a}})}\BibitemShut {NoStop}%
\bibitem [{\citenamefont {Onaga}\ and\ \citenamefont
  {Shinomoto}(2014)}]{PreOnaga}%
  \BibitemOpen
  \bibfield  {author} {\bibinfo {author} {\bibfnamefont {T.}~\bibnamefont
  {Onaga}}\ and\ \bibinfo {author} {\bibfnamefont {S.}~\bibnamefont
  {Shinomoto}},\ }\href {\doibase 10.1103/PhysRevE.89.042817} {\bibfield
  {journal} {\bibinfo  {journal} {Phys. Rev. E}\ }\textbf {\bibinfo {volume}
  {89}},\ \bibinfo {pages} {042817} (\bibinfo {year} {2014})}\BibitemShut
  {NoStop}%
\bibitem [{\citenamefont {Saichev}\ and\ \citenamefont
  {Sornette}(2014)}]{PreSaichev1}%
  \BibitemOpen
  \bibfield  {author} {\bibinfo {author} {\bibfnamefont {A.}~\bibnamefont
  {Saichev}}\ and\ \bibinfo {author} {\bibfnamefont {D.}~\bibnamefont
  {Sornette}},\ }\href {\doibase 10.1103/PhysRevE.89.012104} {\bibfield
  {journal} {\bibinfo  {journal} {Phys. Rev. E}\ }\textbf {\bibinfo {volume}
  {89}},\ \bibinfo {pages} {012104} (\bibinfo {year} {2014})}\BibitemShut
  {NoStop}%
\bibitem [{\citenamefont {Saichev}\ and\ \citenamefont
  {Sornette}(2013)}]{PreSaichev2}%
  \BibitemOpen
  \bibfield  {author} {\bibinfo {author} {\bibfnamefont {A.}~\bibnamefont
  {Saichev}}\ and\ \bibinfo {author} {\bibfnamefont {D.}~\bibnamefont
  {Sornette}},\ }\href {\doibase 10.1103/PhysRevE.87.022815} {\bibfield
  {journal} {\bibinfo  {journal} {Phys. Rev. E}\ }\textbf {\bibinfo {volume}
  {87}},\ \bibinfo {pages} {022815} (\bibinfo {year} {2013})}\BibitemShut
  {NoStop}%
\bibitem [{\citenamefont {Hardiman}\ and\ \citenamefont
  {Bouchaud}(2014)}]{PreBouchaud}%
  \BibitemOpen
  \bibfield  {author} {\bibinfo {author} {\bibfnamefont {S.~J.}\ \bibnamefont
  {Hardiman}}\ and\ \bibinfo {author} {\bibfnamefont {J.-P.}\ \bibnamefont
  {Bouchaud}},\ }\href {\doibase 10.1103/PhysRevE.90.062807} {\bibfield
  {journal} {\bibinfo  {journal} {Phys. Rev. E}\ }\textbf {\bibinfo {volume}
  {90}},\ \bibinfo {pages} {062807} (\bibinfo {year} {2014})}\BibitemShut
  {NoStop}%
\bibitem [{\citenamefont {Dassios}\ and\ \citenamefont
  {Zhao}(2011)}]{DynamicContagion}%
  \BibitemOpen
  \bibfield  {author} {\bibinfo {author} {\bibfnamefont {A.}~\bibnamefont
  {Dassios}}\ and\ \bibinfo {author} {\bibfnamefont {H.}~\bibnamefont {Zhao}},\
  }\href {\doibase 10.1239/aap/1316792671} {\bibfield  {journal} {\bibinfo
  {journal} {Advances in Applied Probability}\ }\textbf {\bibinfo {volume}
  {43}},\ \bibinfo {pages} {814} (\bibinfo {year} {2011})}\BibitemShut
  {NoStop}%
\bibitem [{\citenamefont {Bremaud}\ and\ \citenamefont
  {Massoulie}(1996)}]{Stability}%
  \BibitemOpen
  \bibfield  {author} {\bibinfo {author} {\bibfnamefont {P.}~\bibnamefont
  {Bremaud}}\ and\ \bibinfo {author} {\bibfnamefont {L.}~\bibnamefont
  {Massoulie}},\ }\href@noop {} {\bibfield  {journal} {\bibinfo  {journal} {The
  Annals of Probability}\ }\textbf {\bibinfo {volume} {24}},\ \bibinfo {pages}
  {pp. 1563} (\bibinfo {year} {1996})}\BibitemShut {NoStop}%
\bibitem [{\citenamefont {Zhu}(2012)}]{Zhu2012}%
  \BibitemOpen
  \bibfield  {author} {\bibinfo {author} {\bibfnamefont {L.}~\bibnamefont
  {Zhu}},\ }\href@noop {} {\bibfield  {journal} {\bibinfo  {journal} {Journal
  of Applied Probability}\ ,\ \bibinfo {pages} {760}} (\bibinfo {year}
  {2012})}\BibitemShut {NoStop}%
\bibitem [{\citenamefont {Zhu}(2011)}]{Zhu2011}%
  \BibitemOpen
  \bibfield  {author} {\bibinfo {author} {\bibfnamefont {L.}~\bibnamefont
  {Zhu}},\ }\href@noop {} {\bibfield  {journal} {\bibinfo  {journal} {arXiv
  preprint arXiv:1108.2431}\ ,\ \bibinfo {pages} {1}} (\bibinfo {year}
  {2011})},\ \Eprint {http://arxiv.org/abs/arXiv:1108.2431v2}
  {arXiv:arXiv:1108.2431v2} \BibitemShut {NoStop}%
\bibitem [{\citenamefont {Bartlett}(1963)}]{Bartlett}%
  \BibitemOpen
  \bibfield  {author} {\bibinfo {author} {\bibfnamefont {M.}~\bibnamefont
  {Bartlett}},\ }\href@noop {} {\bibfield  {journal} {\bibinfo  {journal}
  {Journal of the Royal Statistical Society. Series B (Methodological)}\ ,\
  \bibinfo {pages} {264}} (\bibinfo {year} {1963})}\BibitemShut {NoStop}%
\bibitem [{\citenamefont {Hawkes}(1971{\natexlab{b}})}]{Hawkes2}%
  \BibitemOpen
  \bibfield  {author} {\bibinfo {author} {\bibfnamefont {A.~G.}\ \bibnamefont
  {Hawkes}},\ }\href@noop {} {\bibfield  {journal} {\bibinfo  {journal}
  {Journal of the Royal Statistical Society. Series B (Methodological)}\ ,\
  \bibinfo {pages} {438}} (\bibinfo {year} {1971}{\natexlab{b}})}\BibitemShut
  {NoStop}%
\bibitem [{\citenamefont {Adamopoulos}(1975)}]{PGF}%
  \BibitemOpen
  \bibfield  {author} {\bibinfo {author} {\bibfnamefont {L.}~\bibnamefont
  {Adamopoulos}},\ }\href@noop {} {\bibfield  {journal} {\bibinfo  {journal}
  {Journal of Applied Probability}\ }\textbf {\bibinfo {volume} {12}},\
  \bibinfo {pages} {78} (\bibinfo {year} {1975})}\BibitemShut {NoStop}%
\bibitem [{\citenamefont {Saichev}\ \emph {et~al.}(2013)\citenamefont
  {Saichev}, \citenamefont {Maillart},\ and\ \citenamefont
  {Sornette}}]{Saichev2011}%
  \BibitemOpen
  \bibfield  {author} {\bibinfo {author} {\bibfnamefont {A.}~\bibnamefont
  {Saichev}}, \bibinfo {author} {\bibfnamefont {T.}~\bibnamefont {Maillart}}, \
  and\ \bibinfo {author} {\bibfnamefont {D.}~\bibnamefont {Sornette}},\ }\href
  {\doibase 10.1140/epjb/e2013-30493-9} {\bibfield  {journal} {\bibinfo
  {journal} {The European Physical Journal B}\ }\textbf {\bibinfo {volume}
  {86}},\ \bibinfo {eid} {124} (\bibinfo {year} {2013}),\
  10.1140/epjb/e2013-30493-9}\BibitemShut {NoStop}%
\bibitem [{\citenamefont {Saichev}\ and\ \citenamefont
  {Sornette}(2011)}]{Saichev2011a}%
  \BibitemOpen
  \bibfield  {author} {\bibinfo {author} {\bibfnamefont {A.}~\bibnamefont
  {Saichev}}\ and\ \bibinfo {author} {\bibfnamefont {D.}~\bibnamefont
  {Sornette}},\ }\href@noop {} {\bibfield  {journal} {\bibinfo  {journal} {The
  European Physical Journal B-Condensed Matter and Complex Systems}\ }\textbf
  {\bibinfo {volume} {83}},\ \bibinfo {pages} {271} (\bibinfo {year}
  {2011})}\BibitemShut {NoStop}%
\bibitem [{\citenamefont {Rossant}\ \emph {et~al.}(2011)\citenamefont
  {Rossant}, \citenamefont {Leijon}, \citenamefont {Magnusson},\ and\
  \citenamefont {Brette}}]{Rossant2011}%
  \BibitemOpen
  \bibfield  {author} {\bibinfo {author} {\bibfnamefont {C.}~\bibnamefont
  {Rossant}}, \bibinfo {author} {\bibfnamefont {S.}~\bibnamefont {Leijon}},
  \bibinfo {author} {\bibfnamefont {A.~K.}\ \bibnamefont {Magnusson}}, \ and\
  \bibinfo {author} {\bibfnamefont {R.}~\bibnamefont {Brette}},\ }\href@noop {}
  {\bibfield  {journal} {\bibinfo  {journal} {The Journal of Neuroscience}\
  }\textbf {\bibinfo {volume} {31}},\ \bibinfo {pages} {17193} (\bibinfo {year}
  {2011})}\BibitemShut {NoStop}%
\bibitem [{\citenamefont {Kuhn}\ \emph {et~al.}(2003)\citenamefont {Kuhn},
  \citenamefont {Aertsen},\ and\ \citenamefont {Rotter}}]{Kuhn2003}%
  \BibitemOpen
  \bibfield  {author} {\bibinfo {author} {\bibfnamefont {A.}~\bibnamefont
  {Kuhn}}, \bibinfo {author} {\bibfnamefont {A.}~\bibnamefont {Aertsen}}, \
  and\ \bibinfo {author} {\bibfnamefont {S.}~\bibnamefont {Rotter}},\
  }\href@noop {} {\bibfield  {journal} {\bibinfo  {journal} {Neural
  Computation}\ }\textbf {\bibinfo {volume} {15}},\ \bibinfo {pages} {67}
  (\bibinfo {year} {2003})}\BibitemShut {NoStop}%
\bibitem [{\citenamefont {Hawkes}\ and\ \citenamefont {Oakes}(1974)}]{Cluster}%
  \BibitemOpen
  \bibfield  {author} {\bibinfo {author} {\bibfnamefont {A.}~\bibnamefont
  {Hawkes}}\ and\ \bibinfo {author} {\bibfnamefont {D.}~\bibnamefont {Oakes}},\
  }\href@noop {} {\bibfield  {journal} {\bibinfo  {journal} {Journal of Applied
  Probability}\ }\textbf {\bibinfo {volume} {11}},\ \bibinfo {pages} {493}
  (\bibinfo {year} {1974})}\BibitemShut {NoStop}%
\bibitem [{\citenamefont {Cox}\ and\ \citenamefont {Isham}(1980)}]{CoxPoint}%
  \BibitemOpen
  \bibfield  {author} {\bibinfo {author} {\bibfnamefont {D.~R.}\ \bibnamefont
  {Cox}}\ and\ \bibinfo {author} {\bibfnamefont {V.}~\bibnamefont {Isham}},\
  }\href@noop {} {\emph {\bibinfo {title} {Point processes}}},\ Vol.~\bibinfo
  {volume} {12}\ (\bibinfo  {publisher} {CRC Press},\ \bibinfo {year}
  {1980})\BibitemShut {NoStop}%
\bibitem [{\citenamefont {Pernice}\ \emph
  {et~al.}(2011{\natexlab{b}})\citenamefont {Pernice}, \citenamefont {Staude},
  \citenamefont {Cardanobile},\ and\ \citenamefont {Rotter}}]{Volker}%
  \BibitemOpen
  \bibfield  {author} {\bibinfo {author} {\bibfnamefont {V.}~\bibnamefont
  {Pernice}}, \bibinfo {author} {\bibfnamefont {B.}~\bibnamefont {Staude}},
  \bibinfo {author} {\bibfnamefont {S.}~\bibnamefont {Cardanobile}}, \ and\
  \bibinfo {author} {\bibfnamefont {S.}~\bibnamefont {Rotter}},\ }\href@noop {}
  {\bibfield  {journal} {\bibinfo  {journal} {PLoS computational biology}\
  }\textbf {\bibinfo {volume} {7}},\ \bibinfo {pages} {e1002059} (\bibinfo
  {year} {2011}{\natexlab{b}})}\BibitemShut {NoStop}%
\bibitem [{\citenamefont {Rasmussen}(2013)}]{Bayes}%
  \BibitemOpen
  \bibfield  {author} {\bibinfo {author} {\bibfnamefont {J.}~\bibnamefont
  {Rasmussen}},\ }\href {\doibase 10.1007/s11009-011-9272-5} {\bibfield
  {journal} {\bibinfo  {journal} {Methodology and Computing in Applied
  Probability}\ }\textbf {\bibinfo {volume} {15}},\ \bibinfo {pages} {623}
  (\bibinfo {year} {2013})}\BibitemShut {NoStop}%
\bibitem [{\citenamefont {Lukacs}(1970)}]{lukacs}%
  \BibitemOpen
  \bibfield  {author} {\bibinfo {author} {\bibfnamefont {E.}~\bibnamefont
  {Lukacs}},\ }\href@noop {} {\emph {\bibinfo {title} {Characteristic
  Functions}}},\ Griffin books of Cognate Interest\ (\bibinfo  {publisher}
  {Hafner Publishing Company},\ \bibinfo {year} {1970})\BibitemShut {NoStop}%
\bibitem [{\citenamefont {Bacry}\ \emph {et~al.}(2012)\citenamefont {Bacry},
  \citenamefont {Dayri},\ and\ \citenamefont {Muzy}}]{bacry2012}%
  \BibitemOpen
  \bibfield  {author} {\bibinfo {author} {\bibfnamefont {E.}~\bibnamefont
  {Bacry}}, \bibinfo {author} {\bibfnamefont {K.}~\bibnamefont {Dayri}}, \ and\
  \bibinfo {author} {\bibfnamefont {J.-F.}\ \bibnamefont {Muzy}},\ }\href@noop
  {} {\bibfield  {journal} {\bibinfo  {journal} {The European Physical Journal
  B-Condensed Matter and Complex Systems}\ }\textbf {\bibinfo {volume} {85}},\
  \bibinfo {pages} {1} (\bibinfo {year} {2012})}\BibitemShut {NoStop}%
\bibitem [{\citenamefont {Felsenstein}(2004)}]{Felsenstein}%
  \BibitemOpen
  \bibfield  {author} {\bibinfo {author} {\bibfnamefont {J.}~\bibnamefont
  {Felsenstein}},\ }\href@noop {} {\emph {\bibinfo {title} {Inferring
  Phylogenies}}}\ (\bibinfo  {publisher} {Sinauer Associates, Incorporated},\
  \bibinfo {year} {2004})\BibitemShut {NoStop}%
\bibitem [{\citenamefont {Hebb}(1949)}]{hebb1949}%
  \BibitemOpen
  \bibfield  {author} {\bibinfo {author} {\bibfnamefont {D.~O.}\ \bibnamefont
  {Hebb}},\ }\href@noop {} {\emph {\bibinfo {title} {{The Organization of
  Behavior: A Neuropsychological Theory}}}},\ \bibinfo {edition} {new ed}\ ed.\
  (\bibinfo  {publisher} {Wiley},\ \bibinfo {address} {New York},\ \bibinfo
  {year} {1949})\BibitemShut {NoStop}%
\bibitem [{\citenamefont {Gerstein}\ \emph {et~al.}(1989)\citenamefont
  {Gerstein}, \citenamefont {Bedenbaugh},\ and\ \citenamefont
  {Aertsen}}]{gerstein1989}%
  \BibitemOpen
  \bibfield  {author} {\bibinfo {author} {\bibfnamefont {G.}~\bibnamefont
  {Gerstein}}, \bibinfo {author} {\bibfnamefont {P.}~\bibnamefont
  {Bedenbaugh}}, \ and\ \bibinfo {author} {\bibfnamefont {A.~M.}\ \bibnamefont
  {Aertsen}},\ }\href@noop {} {\bibfield  {journal} {\bibinfo  {journal}
  {Biomedical Engineering, IEEE Transactions on}\ }\textbf {\bibinfo {volume}
  {36}},\ \bibinfo {pages} {4} (\bibinfo {year} {1989})}\BibitemShut {NoStop}%
\bibitem [{\citenamefont {Gray}\ and\ \citenamefont {Singer}(1989)}]{Gray1989}%
  \BibitemOpen
  \bibfield  {author} {\bibinfo {author} {\bibfnamefont {C.~M.}\ \bibnamefont
  {Gray}}\ and\ \bibinfo {author} {\bibfnamefont {W.}~\bibnamefont {Singer}},\
  }\href {\doibase 10.1073/pnas.86.5.1698} {\bibfield  {journal} {\bibinfo
  {journal} {Proceedings of the National Academy of Sciences}\ }\textbf
  {\bibinfo {volume} {86}},\ \bibinfo {pages} {1698} (\bibinfo {year}
  {1989})}\BibitemShut {NoStop}%
\bibitem [{\citenamefont {Vaadia}\ \emph {et~al.}(1995)\citenamefont {Vaadia},
  \citenamefont {Haalman}, \citenamefont {Abeles}, \citenamefont {Bergman},
  \citenamefont {Prut}, \citenamefont {Slovin},\ and\ \citenamefont
  {Aertsen}}]{vaadia1995}%
  \BibitemOpen
  \bibfield  {author} {\bibinfo {author} {\bibfnamefont {E.}~\bibnamefont
  {Vaadia}}, \bibinfo {author} {\bibfnamefont {I.}~\bibnamefont {Haalman}},
  \bibinfo {author} {\bibfnamefont {M.}~\bibnamefont {Abeles}}, \bibinfo
  {author} {\bibfnamefont {H.}~\bibnamefont {Bergman}}, \bibinfo {author}
  {\bibfnamefont {Y.}~\bibnamefont {Prut}}, \bibinfo {author} {\bibfnamefont
  {H.}~\bibnamefont {Slovin}}, \ and\ \bibinfo {author} {\bibfnamefont
  {A.}~\bibnamefont {Aertsen}},\ }\href@noop {} {\bibfield  {journal} {\bibinfo
   {journal} {Nature}\ }\textbf {\bibinfo {volume} {373}},\ \bibinfo {pages}
  {515} (\bibinfo {year} {1995})}\BibitemShut {NoStop}%
\bibitem [{\citenamefont {Riehle}\ \emph {et~al.}(1997)\citenamefont {Riehle},
  \citenamefont {Gr{\"u}n}, \citenamefont {Diesmann},\ and\ \citenamefont
  {Aertsen}}]{riehle1997}%
  \BibitemOpen
  \bibfield  {author} {\bibinfo {author} {\bibfnamefont {A.}~\bibnamefont
  {Riehle}}, \bibinfo {author} {\bibfnamefont {S.}~\bibnamefont {Gr{\"u}n}},
  \bibinfo {author} {\bibfnamefont {M.}~\bibnamefont {Diesmann}}, \ and\
  \bibinfo {author} {\bibfnamefont {A.}~\bibnamefont {Aertsen}},\ }\href@noop
  {} {\bibfield  {journal} {\bibinfo  {journal} {Science}\ }\textbf {\bibinfo
  {volume} {278}},\ \bibinfo {pages} {1950} (\bibinfo {year}
  {1997})}\BibitemShut {NoStop}%
\bibitem [{\citenamefont {Bair}\ \emph {et~al.}(2001)\citenamefont {Bair},
  \citenamefont {Zohary},\ and\ \citenamefont {Newsome}}]{bair2001}%
  \BibitemOpen
  \bibfield  {author} {\bibinfo {author} {\bibfnamefont {W.}~\bibnamefont
  {Bair}}, \bibinfo {author} {\bibfnamefont {E.}~\bibnamefont {Zohary}}, \ and\
  \bibinfo {author} {\bibfnamefont {W.~T.}\ \bibnamefont {Newsome}},\
  }\href@noop {} {\bibfield  {journal} {\bibinfo  {journal} {The journal of
  Neuroscience}\ }\textbf {\bibinfo {volume} {21}},\ \bibinfo {pages} {1676}
  (\bibinfo {year} {2001})}\BibitemShut {NoStop}%
\bibitem [{\citenamefont {Kohn}\ and\ \citenamefont {Smith}(2005)}]{kohn2005}%
  \BibitemOpen
  \bibfield  {author} {\bibinfo {author} {\bibfnamefont {A.}~\bibnamefont
  {Kohn}}\ and\ \bibinfo {author} {\bibfnamefont {M.~A.}\ \bibnamefont
  {Smith}},\ }\href@noop {} {\bibfield  {journal} {\bibinfo  {journal} {The
  Journal of neuroscience}\ }\textbf {\bibinfo {volume} {25}},\ \bibinfo
  {pages} {3661} (\bibinfo {year} {2005})}\BibitemShut {NoStop}%
\bibitem [{\citenamefont {Martignon}\ \emph {et~al.}(1995)\citenamefont
  {Martignon}, \citenamefont {Von~Hassein}, \citenamefont {Gr{\"u}n},
  \citenamefont {Aertsen},\ and\ \citenamefont {Palm}}]{martignon1995}%
  \BibitemOpen
  \bibfield  {author} {\bibinfo {author} {\bibfnamefont {L.}~\bibnamefont
  {Martignon}}, \bibinfo {author} {\bibfnamefont {H.}~\bibnamefont
  {Von~Hassein}}, \bibinfo {author} {\bibfnamefont {S.}~\bibnamefont
  {Gr{\"u}n}}, \bibinfo {author} {\bibfnamefont {A.}~\bibnamefont {Aertsen}}, \
  and\ \bibinfo {author} {\bibfnamefont {G.}~\bibnamefont {Palm}},\ }\href@noop
  {} {\bibfield  {journal} {\bibinfo  {journal} {Biological cybernetics}\
  }\textbf {\bibinfo {volume} {73}},\ \bibinfo {pages} {69} (\bibinfo {year}
  {1995})}\BibitemShut {NoStop}%
\bibitem [{\citenamefont {Boht{\'e}}\ \emph {et~al.}(2000)\citenamefont
  {Boht{\'e}}, \citenamefont {Spekreijse},\ and\ \citenamefont
  {Roelfsema}}]{bohte2000}%
  \BibitemOpen
  \bibfield  {author} {\bibinfo {author} {\bibfnamefont {S.~M.}\ \bibnamefont
  {Boht{\'e}}}, \bibinfo {author} {\bibfnamefont {H.}~\bibnamefont
  {Spekreijse}}, \ and\ \bibinfo {author} {\bibfnamefont {P.~R.}\ \bibnamefont
  {Roelfsema}},\ }\href@noop {} {\bibfield  {journal} {\bibinfo  {journal}
  {Neural Computation}\ }\textbf {\bibinfo {volume} {12}},\ \bibinfo {pages}
  {153} (\bibinfo {year} {2000})}\BibitemShut {NoStop}%
\bibitem [{\citenamefont {Ohiorhenuan}\ \emph {et~al.}(2010)\citenamefont
  {Ohiorhenuan}, \citenamefont {Mechler}, \citenamefont {Purpura},
  \citenamefont {Schmid}, \citenamefont {Hu},\ and\ \citenamefont
  {Victor}}]{Nature2010}%
  \BibitemOpen
  \bibfield  {author} {\bibinfo {author} {\bibfnamefont {I.~E.}\ \bibnamefont
  {Ohiorhenuan}}, \bibinfo {author} {\bibfnamefont {F.}~\bibnamefont
  {Mechler}}, \bibinfo {author} {\bibfnamefont {K.~P.}\ \bibnamefont
  {Purpura}}, \bibinfo {author} {\bibfnamefont {A.~M.}\ \bibnamefont {Schmid}},
  \bibinfo {author} {\bibfnamefont {Q.}~\bibnamefont {Hu}}, \ and\ \bibinfo
  {author} {\bibfnamefont {J.~D.}\ \bibnamefont {Victor}},\ }\href@noop {}
  {\bibfield  {journal} {\bibinfo  {journal} {Nature}\ }\textbf {\bibinfo
  {volume} {466}},\ \bibinfo {pages} {617} (\bibinfo {year}
  {2010})}\BibitemShut {NoStop}%
\bibitem [{\citenamefont {Cardanobile}\ and\ \citenamefont
  {Rotter}(2009)}]{Cardanobile2009}%
  \BibitemOpen
  \bibfield  {author} {\bibinfo {author} {\bibfnamefont {S.}~\bibnamefont
  {Cardanobile}}\ and\ \bibinfo {author} {\bibfnamefont {S.}~\bibnamefont
  {Rotter}},\ }\href@noop {} {\bibfield  {journal} {\bibinfo  {journal}
  {Journal of Computational Neuroscience}\ }\textbf {\bibinfo {volume} {28}},\
  \bibinfo {pages} {267} (\bibinfo {year} {2009})},\ \Eprint
  {http://arxiv.org/abs/arXiv:0904.1505v3} {arXiv:0904.1505v3} \BibitemShut
  {NoStop}%
\end{thebibliography}
%merlin.mbs apsrev4-1.bst 2010-07-25 4.21a (PWD, AO, DPC) hacked
%Control: key (0)
%Control: author (8) initials jnrlst
%Control: editor formatted (1) identically to author
%Control: production of article title (-1) disabled
%Control: page (0) single
%Control: year (1) truncated
%Control: production of eprint (0) enabled
%

\end{document}